\documentclass[reqno,a4paper,11pt]{amsart}

\usepackage[all]{xy}
\usepackage{amsfonts}
\usepackage[mathcal]{eucal}
\usepackage{eufrak}
\usepackage{amssymb}
\usepackage{amsmath}
\usepackage{color}
\usepackage[pagebackref,
colorlinks
]{hyperref}

\newcommand{\ov}{\overline}

\newcommand{\dlambda}{\delta}
\newcommand{\id}{\operatorname{id}}

\newcommand{\nn}{\mathbb N}

\newcommand{\sE}{\mathsf{E}}
\newcommand{\sT}{\mathsf{T}}

\newcommand{\sB}{\mathsf{B}}
\newcommand{\sC}{\mathsf{C}}
\newcommand{\sL}{\mathsf{L}}
\newcommand{\sM}{\mathsf{M}}

\newcommand{\sA}{\mathsf{S}}

\newcommand{\Ga}{\Gamma}
\newcommand{\ga}{\gamma}

\newcommand{\de}{\delta}
\newcommand{\la}{\lambda}

\newcommand\SK[1][1]{\operatorname{SK}_{#1}}
\newcommand\hSK[1][1]{\operatorname{SK}^{h}_{#1}}
\newcommand{\Nrd}[1][{}]{{\operatorname{Nrd}_{#1}}}
\newcommand{\Trd}[1][{}]{{\operatorname{Trd}_{#1}}}

\newcommand{\GL}{\operatorname{GL}}

\newcommand{\sgn}{\operatorname{sgn}}
\newcommand{\M}{\mathbb M}
\newcommand{\I}{\mathbb I}
\newcommand{\D}{\mathbb D}

\newcommand{\chr}{\operatorname{char}}
\newcommand{\ind}{\operatorname{ind}}

\newcommand{\Gal}{\operatorname{Gal}}

\newcommand{\gr}{\text{gr}}

\newcommand{\inv}{^{-1}}
\newcommand{\diag}{\operatorname{diag}}
\newcommand{\elh}{\mathcal E \ell_h}

\newcommand{\gbr}{\operatorname{grBr}}
\newcommand{\br}{\operatorname{Br}}

\theoremstyle{plain}
\newtheorem {lemma}{Lemma}[section]
\newtheorem {theorem}[lemma]{Theorem}

\newtheorem {cor}[lemma]{Corollary}
\newtheorem {proposition}[lemma]{Proposition}

\theoremstyle{remark}
\newtheorem* {remark}{Remark}  % for unnumbered remark
\newtheorem {example}[lemma]{Example}

\theoremstyle{definition}

\numberwithin{equation}{section}

\title[Homogeneous $\SK$ of simple graded algebras]{Homogeneous $\mathbf {SK_1}$ of simple graded algebras}

\date{}

\begin{thanks}
{This work was supported by the Engineering and Physical
  Sciences Research Council [grant EP/I007784/1].}
\end{thanks}

\author{R. Hazrat}\address{
Department of Pure Mathematics\\
Queen's University\\
Belfast BT7 1NN\\
United Kingdom} 
\email{r.hazrat@qub.ac.uk}

\author{A. R. Wadsworth}
\address{
Department of Mathematics\\
University of California at San Diego\\
La Jolla, California 92093-0112\\
U.S.A.}
 \email{arwadsworth@ucsd.edu}

\begin{document}

\begin{abstract}
For a simple graded algebra $\sA=\M_n(\sE)$ over a 
graded division algebra $\sE$, a short exact sequence 
relating the reduced Whitehead group of the homogeneous 
part of $\sA$ to that of $\sE$ is established. In 
particular it is shown that the homogeneous $\SK$ is 
not in general  Morita invariant. 
\end{abstract}

\maketitle

Graded methods in the theory of valued division algebras 
have  proved to be extremely useful. A valuation $v$
 on a division algebra $D$ induces a filtration 
on $D$ which yields an associated graded ring~
$\gr(D)$.  Indeed, $\gr(D)$ is a graded division algebra,
i.e., every nonzero homogeneous element of $\gr(D)$
is a unit. While $\gr(D)$ has a much simpler
structure than $D$, nonetheless $\gr(D)$ provides a 
remarkably good  reflection of $D$ in many ways, 
particualrly when the valuation on the center 
$Z(D)$ is Henselian.  The approach of making calculations
in $\gr(D)$, then lifting back to get nontrivial 
information about $D$ has been remarkably successful.
See \cite{jw, wadval} for background on valued division 
algebras,
and \cite{hwcor,tigwad1,tigwad2} for connections  between 
valued and graded division algebras. The recent papers
\cite{HW1,HW2,WY,WSK1} on the reduced Whitehead group $\SK$ and its
unitary analogue  have provided good illustrations of 
the effectiveness of this approach.  Notably it was proved 
in \cite[Th.~4.8, Th.~5.7]{HW1}  that if $v$ on $Z(D)$ is 
Henselian and $D$ is tame over $Z(D)$, then $\SK(D) \cong
\SK(\gr(D))$ and $\SK(\gr(D))\cong \SK(q(\gr(D)))$, 
where $q(\gr(D))$ is the division ring of quotients of 
$\gr(D)$.  This has allowed recovery of many of the known 
calculations of $\SK(D)$ with much easier proofs, as well
as leading to determinations of $\SK(D)$ in some
new cases.

By the graded Wedderburn theorem 
(see~\cite[Prop.~1.3(a)]{hwcor}), 
any simple graded algebra $\sA$
 finite-dimensional over its center $\sT$  has the form 
$\sA=\M_n (\sE)(\de_1 , \ldots , \de_n)$, 
where the $\delta_i$ lie in
 a torsion-free abelian group $\Gamma$ containing the 
grade group $\Gamma_\sE$.  That is $\sA$ is the 
 $n \times n$ matrix algebra over 
a graded division algebra~$\sE$ with its grading shifted 
by $(\de_1 , \ldots , \de_n)$. Since $\sA$ is known  to 
be   Azumaya 
algebra over $\sT$, there is  a reduced norm  
map on the group of units, $\Nrd_\sA\colon \sA^* \to  
\sT^*$;  one can then  
define the reduced Whitehead group 
$\SK(\sA)$ in the usual manner as the kernel of the 
reduced norm of $\sA$ modulo the commutator subgroup of 
$\sA^*$ (see Def.~\ref{azuski}). However $\SK$ is not 
a ``graded functor'', i.e., it does not take into account 
the grading on~$\sA$. 

To factor in   the grading on  $\sA$,  we introduce 
in this paper the {\it homogeneous} reduced
 Whitehead group $\hSK(\sA)$ (see Def.~\ref{zuhgoi}),
which treats only the homogeneous units of $\sA$. 
We
establish a short exact sequence relating   
$\hSK(\sA)$ to $\SK(\sE)$ (see Th.~\ref{extseq}) which   
 allows us to calculate $\hSK(\sA)$ in many cases.  
In particular we show that $\hSK$ is not in general 
Morita invariant for $\sE$, and indeed can behave quite badly 
when the semisimple ring $\sA_0$ is not simple
(see Ex.~\ref{hgdstr}). 
As a prelude to this, in ~\S 1 we prove the existence and
multiplicativity of a Dieudonn\'e determinant for homogeneous
elements of $\sA = \M_n(\sE)(\delta_1, \ldots, \de_n)$.  
This was originally 
needed for the work on $\SK^h$, but later it turned out that 
the ungraded Dieutdonn\'e determinant for the semisimple algebra
$\sA_0$ was all that was needed.  We have nonetheless included 
the development of the homogoneous Dieudonn\'e determinant, since we 
feel that it is of some interest in its own right.

\section{Dieudonn\'{e} determinant}

Throughout this paper we will be working with 
matrices over graded 
division rings. Recall that a graded ring
$\sE = \bigoplus_{\gamma\in \Gamma} \sE_\gamma$ 
is called a {\it graded division ring} if every 
non-zero homogeneous element of $\sE$ is a unit, i.e. it has a (two-sided)
multiplicative inverse. We assume throughout that
the index set $\Gamma$ is an abelian group.  
Note that the hypothesis on $\sE$ implies that  the 
grade set $\Ga_{\sE} =\{\gamma\in \Gamma\mid
\sE_\gamma \ne \{0\}\}$ is actually a subgroup of $\Gamma$.
We write $\sE_h^*$ for the group of homogeneous units
of $\sE$, which consists of all the nonzero homogeneous
elements of $\sE$, and can be a proper subgroup of the 
group $\sE^*$ of all units of $\sE$. 
  
 Let $\M_n(\sE)$ be the $n\times n$ matrix ring over
the graded division ring $\sE$.
For any $x\in \sE$, let
$E_{ij}(x)$ be the matrix in $\M_n(\sE)$ with  $x$ in 
$(i,j)$-position and $0$'s otherwise.
Take any $\de_1, \ldots, \de_n \in \Gamma$.
  The shifted 
grading on $\M_n(\sE)$ determined by 
$(\de_1, \ldots, \de_n)$ is defined by setting, 
\begin{equation}\label{hogr}
\deg(E_{ij}(x))=\deg(x)+\delta_i-\delta_j\qquad
\text{for any homogeneous $x$ in $\sE$}.
\end{equation}  
This is then extended linearly to all of 
$\M_n(\sE)$.
 One can then see that  for $\la \in \Gamma$, 
the $\la$-component ${\M_n (\sE)}_{\la}$  consists of
those matrices with homogeneous entries, with the 
degrees shifted as follows:
\begin{equation}\label{mmkkhh}
{\M_n(\sE)}_{\la}\, =\,
\begin{pmatrix}
\sE_{ \la+\de_1 - \de_1} & \sE_{\la+\de_2  - \de_1} & \cdots &
\sE_{\la +\de_n - \de_1} \\
\sE_{\la + \de_1 - \de_2} & \sE_{\la + \de_2 - \de_2} & \cdots &
\sE_{\la+\de_n  - \de_2} \\
\vdots  & \vdots  & \ddots & \vdots  \\
\sE_{\la + \de_1 - \de_n} & \sE_{ \la + \de_2 - \de_n} & \cdots &
\sE_{\la + \de_n - \de_n}
\end{pmatrix}.
\end{equation}
That is, ${\M_n(\sE)}_{\la}$ consists of matrices with 
each  $ij$-entry lying in
$\sE_{\la+\delta_j  - \delta_i}$. 
We then have 
\[
\textstyle
\M_n (\sE)\,=\,\bigoplus\limits_{\la \in \Gamma} 
{\M_n (\sE)}_{\la} \quad \text {and}
\quad \M_n (\sE)_\lambda \cdot \M_n (\sE)_\mu
\,\subseteq\, \M_n (\sE)_{\lambda +\mu}
\ \ \text{for all}\ \ \lambda,\mu \in \Gamma,
\] 
which shows that $\M_n(\sE)$ is a graded ring.
We denote the matrix ring with this grading by 
$\M_n (\sE)(\de_1 , \ldots , \de_n)$ or 
$\M_n (\sE)(\overline \de)$, where 
$\overline \de=(\de_1 , \ldots , \de_n)$.
It is not hard to show that $\M_n(\sE)(\ov \de)$
is a simple graded ring, i.e., it has no nontrivial
homogeneous two-sided ideals.  
Observe that for $\sA = \M_n(\sE(\ov \de))$, the grade
set is
\begin{equation}\label{GammaS}
\textstyle 
\Gamma_\sA \, = \, \bigcup \limits_{i= 1}^n 
\bigcup\limits
_{j= 1}^n \,(\delta_j - \delta _i) + \Gamma_\sE,
\end{equation}
which need not be a group.  However, if we let 
$\sA_h^* = \{ A \in \sA \mid A \text{ is homogeneous  and
 $A$ is a unit of $\sA$}\}$, which is a subgroup of 
the group of units $\sA^*$ of $\sA$, and set
$$
\Gamma_\sA^* \, = \,\{\deg(A) \mid A \in \sA_h^*\},
$$
then $\Gamma_\sA^*$ is a subgroup of $\Gamma$, with 
$\Gamma_\sE \subseteq \Gamma_\sA^* \subseteq \Ga_\sA$.  

Note that 
when $\de_i=0$, $1\leq i \leq n$, then 
${\M_n(\sE)}_{\la}={\M_n(\sE_{\la})}$. We refer to this 
case as the {\it unshifted grading} on $\M_n(\sE)$.

For any graded rings $\sB$ and $\sC$, we write 
$\sB \cong_\gr \sC$ if there is graded ring isomorphism
$\sB \to \sC$, i.e., a ring isomorphism 
 that  maps
$\sB_\lambda$ onto $\sC_\lambda$ for all $\lambda \in 
\Gamma_\sB = \Gamma_\sC$.

The following two statements can be proved easily 
(see~\cite[pp.~60-61]{grrings}):
\begin{itemize}
\item[$\circ$] If $\alpha \in \Gamma$, and 
$\pi \in S_n$ is a permutation then 
\begin{equation} \label{pqow1}
\M_n (\sE)(\de_1 , \ldots , \de_n)\,\cong_{\gr} \,
\M_n (\sE)(\de_{\pi(1)}+\alpha , \ldots , 
\de_{\pi(n)}+\alpha).
\end{equation}

\item[$\circ$]   If $\alpha_1,\dots, \alpha_n \in \Ga$ 
with $\alpha_i=\deg(u_i)$ for some units $u_i \in \sE_h^*$, 
then 
\begin{equation}\label{pqow2}
\M_n (\sE)(\de_1 , \ldots , \de_n)\,\cong_{\gr}\,
 \M_n (\sE)(\de_1+\alpha_1 , \ldots , \de_n+\alpha_n).
\end{equation}
\end{itemize}

Take any $\de_1, \ldots, \de_n\in \Gamma$ 
In the factor group $\Ga/\Gamma_{\sE}$ and 
let $\varepsilon_1 +\Ga_{\sE},\dots,\varepsilon_k+
\Ga_{\sE}$  be the distinct cosets  in 
 $\{\de_1+\Ga, \ldots, \de_n+ \Ga\}$. 
 For each $\varepsilon_\ell$, 
let $r_\ell$ 
be the number of $i$ with $\de_i+\Ga=\varepsilon_\ell+\Ga$.  
It was observed in \cite[Prop.~1.4]{hwcor} that 
\begin{equation}\label{urnha}
\M_n(\sE)_0 \,\cong\, \M_{r_1}(\sE_0)\times \ldots 
\times \M_{r_k}(\sE_0).
\end{equation}
Thus  $\M_n(\sE)_0$ 
is a a semisimple  ring; it is simple if and 
only if $k=1$. Indeed, (\ref{urnha}) follows easily 
from the observations above.  For, 
using (\ref{pqow1}) and 
(\ref{pqow2}) we get 
\begin{equation}\label{wdeild}
\M_n(\sE)(\de_1,\dots,\de_n)\,\cong_{\gr}\,\M_n(\sE)
(\varepsilon_1,\dots,\varepsilon_1,\varepsilon_2,\dots,
\varepsilon_2,\dots,\varepsilon_k,\dots,\varepsilon_k),
\end{equation}
with each $\varepsilon_\ell$ occurring $r_\ell$ times.  
Now (\ref{mmkkhh}) for $\lambda=0$ and 
$(\delta_1,\dots,\delta_n)=
(\varepsilon_1,\dots\varepsilon_1,\varepsilon_2,\dots,
\varepsilon_2,\dots,\varepsilon_k,\dots,\varepsilon_k)$ 
immediately gives (\ref{urnha}).

If the graded ring $\sE$ is commutative then the usual determinant 
map is available, and ${\det\big({\M_n(\sE)}_{\la}\big) 
\subseteq \sE_{n\la}}$. Indeed, if 
$a=(a_{ij}) \in {\M_n(\sE)}_{\la}$, then 
$\det(a)=
\sum_{\sigma \in S_n}\sgn(\sigma)a_{1\sigma1}a_{2\sigma2}
\dots a_{n\sigma n}\in \sE$. But by ~(\ref{mmkkhh}) 
\begin{equation}\label{jsjgfd}
\textstyle
\deg(a_{1\sigma1}a_{2\sigma2}\dots a_{n\sigma n})\, =\ 
n\la + \sum\limits_{i=1}^n \de_{\sigma(i)} - 
\sum\limits_{i=1}^n \de_i\ = \ n\la.
\end{equation} 
When  $\sE$  is not commutative, there is no 
well-defined determinant available in general. For a 
division ring $D$, Dieudonn\'{e} constructed a 
determinant map which reduces to the usual 
  determinant when $D$ is commutative. This is a 
group homomorphism 
$\det\colon\GL_n(D) \rightarrow D^*/[D^*,D^*]$.
 The kernel of $\det$ is  the 
subgroup $E_n(D)$ of $\GL_n(D)$ generated by 
elementary matrices,
 which coincides with 
the commutator group $[\GL_n(D),\GL_n(D)]$ 
unless $\M_n(D) = \M_2(\mathbb F_2)$ 
(see Draxl~\cite[\S20]{draxl}). Note that the 
construction of a Dieudonn\'e
determinant has  been carried over to 
(noncommutative) local and semilocal rings 
in (\cite{vas}). 

Since graded division rings behave in many ways like local
rings, one may ask whether there is a map like 
the Dieudonn\'e determinant in the graded setting.
We will show that this is indeed the case,
so long as one restricts to homogeneous elements.
Specifically, let $\sE$ be a graded division ring with 
 grade group 
$\Gamma_E\subseteq \Gamma$ with $\Gamma$ abelian, and let
$\overline \de=(\de_1,\dots, \de_n)$, where 
$\de_i \in \Gamma$. Let $\sA=\M_n(\sE)(\overline \de)$ 
be the matrix ring over $\sE$ with 
grading shifted by $\overline \de$.  
Denote by $\sA_h$ the set of homogeneous 
elements of $\sA$ and by $\sA^*_h$ or 
$\GL_n^h(\sE)(\overline \de)$ the 
group of  homogeneous units of $\sA$. 
We will show in Th.~\ref{homde} that there is a determinant-like  group 
homomorphism  $\det_\sE\colon\sA^*_h \rightarrow 
\sE_h^*/[\sE_h^*,\sE_h^*]$ 
which is compatible with the Dieudonn\'{e} determinant 
  on the semisimple 
ring $\sA_0$ (see commutative diagram \eqref{ddetE0}).

We first show that every matrix in 
$\GL_n^h(\sE)(\overline \de)$ can be  decomposed into  
strict  Bruhat normal form.  In this decomposition, a
triangular matrix  is said to be {\it unipotent}
triangular if all its diagonal entries are $1$'s.

\begin{proposition}\label{bruhatf} $($Bruhat normal form$)$ 
Let $\sE$ be a graded division ring with  grade group 
$\Gamma_\sE \subseteq\Gamma$. 
Let ${\sA=\M_n(\sE)(\overline \de)}$ be a 
matrix ring over $\sE$ with grading shifted by 
$\overline \de=(\de_1,\dots, \de_n)$,  $\de_i \in \Gamma$. Then
every  $A\in  \sA_h^*$ has a unique strict Bruhat 
normal form, 
i.e., $A$ can be decomposed uniquely as 
$$
A\,=\,TUP_{\pi}V
$$
 for matrices
$T,U, P_\pi, V$ in $\sA$ such that $T$ is unipotent lower
 triangular,
$U$ is diagonal and invertible, $P_{\pi}$ is a permutation  
matrix, and $V$ is unipotent
upper triangular with $P_{\pi}VP_{\pi}^{-1}$ also 
unipotent upper triangular.
Moreover, $T$, $UP_{\pi}$, and $V$ are homogeneous 
matrices, with 
$\deg(T)=  \deg(V)= 0$ and $\deg(UP_{\pi})= \deg(A)$.
Also, $T$ is a product of homogeneous elementary matrices
$($of degree $0)$.
\end{proposition}
 
\begin{proof}
The construction follows closely that in  Draxl 
\cite[\S19, Thm~1]{draxl},
 with extra attention given to degrees of the homogeneous 
matrices in 
the graded ring $\sA = \M_n(\sE)(\overline \de)$.
We will  carry out elementary row operations on 
homogeneous invertible matrices, which corresponds to 
left multiplication by elementary matrices.  But, we use 
only homogeneous elementary matrices thereby preserving 
homogeneity of the matrices being reduced. 
For $x\in \sE$ and $i,j \in \{1, 2, \ldots, n\}$ with 
$i\ne j$, let $e_{ij}(x) = \I_n + E_{ij}(x)$, which is 
the elementary matrix with all $1$'s on the main diagonal,
$x$ in the $(i,j)$-position and all other entries $0$.
Note that if $e_{ij}(x)$ is homogeneous, it must have 
degree $0$ because of the $1$'s on the main diagonal.
So, in view of \eqref{mmkkhh}, $e_{ij}(x)$ is homogeneous
if and only if $x$ is homogeneous with $\deg(x) = 
\delta_j - \de _i$ or $x = 0$.  
Let
\begin{equation}\label{elh}
\elh\, = \, \{\text{homogeneous elementary matrices in
 $\sA$}\} \, = \,
\{ e_{ij}(x)\mid i\ne j \text{ and } x\in 
\sE_{\delta_j-\delta_i}\}.
\end{equation}

Let $A\in \sA_h^*$. Since $A$ is homogeneous, 
every  non-zero 
entry of $A$ is a  homogeneous element of the graded 
division ring 
$\sE$ (see~\eqref{mmkkhh}), and so is  a unit of $\sE$. 
Since $A$ is an invertible matrix, each row must have at 
least one 
nonzero entry.  
Write the $(i,j)$-entry of $A$ as $a^1_{ij}$; so 
$A=(a^1_{ij})$.
  Let $a^1_{1\rho(1)}$ be the first nonzero entry in the 
first row, 
working from the left. For $i >1$, multiplying 
$A$ on the left by the elementary matrix  
$e_{i1}(-a^1_{i\rho(1)}(a^1_{1\rho(1)})^{-1})$ amounts to adding  
the left multiple $-a^1_{i\rho(1)}(a^1_{1\rho(1)})^{-1}$ times the 
first  row to the $i$-th row;  it makes the $(i,\rho(1))$-entry 
zero, without altering any other rows besides the $i$-th.  
%% Repeating this process $n-1$ times, for $2\leq i \leq n$, we can 
%% make all the entries below $a_{1\rho(1)}$ zero. 
By iterating this for each row below the first row, we 
obtain a matrix $A^{(1)}=\prod_{i=n}^2 
e_{i1}(-a^1_{i\rho(1)}(a^1_{1\rho(1)})^{-1})A$, 
which  has the form 
\begin{equation}\label{moirt}
A^{(1)}=
\begin{pmatrix}
0 & 0                   & \cdots & a^1_{1\rho(1)} & a^1_{1,\rho(1)+1}& \cdots & a^1_{1n} \\
a^1_{21} & a^1_{22}  & \cdots &      0             & b_{2,\rho(1)+1}& \cdots & b_{2n} \\
a^1_{31} & a^1_{32}  & \cdots &      0             & b_{3,\rho(1)+1}& \cdots & b_{3n} \\
\vdots  & \vdots   & \vdots &      \vdots            &\vdots               & \vdots  & \vdots \\
a^1_{n1} & a^1_{n2}  & \cdots &      0             & b_{n,\rho(1)+1}& \cdots & b_{nn} \\
\end{pmatrix}.
\end{equation}
Let  $\la = \deg(A)$. From the definition of the  grading on 
$\M_n(A)(\overline \delta)$  
we have   ${\deg(a^1_{i\rho(1)})=\la+\de_{\rho(1)}-\de_i}$ (see~\ref{mmkkhh}). 
Therefore $\deg(-a^1_{i\rho(1)}(a^1_{1\rho(1)})^{-1})
=\de_1-\de_i$, 
which shows that  $e_{i1}(-a^1_{i\rho(1)}(a^1_{1\rho(1)})^{-1})\in
\elh$ for $i = 2,3, \ldots, n$.  
Since homogeneous elementary matrices have degree $0$, 
$A^{(1)}$ is homogeneous with 
$\deg(A^{(1)})=\deg(A)=\la$. 

Write $A^{(1)}=(a^2_{ij})$. Since $A^{(1)}$ is 
invertible, not all the entries of its second row can be  
zero. Let $a^2_{2\rho(2)}$ be the first nonzero entry 
 in the second row working from the left
(clearly $\rho(1) \not = \rho(2)$), and repeat 
the  process above  with $A^{(1)}$ to get a 
homogeneous invertible matrix 
$A^{(2)}$ with all entries below $a^2_{2\rho(2)}$ zero.
In doing this, the entries in the $\rho(1)$ column are
unchanged. 
By iterating this process, working down 
row by row, we obtain
a matrix 
$A^{(n-1)}=\big(a_{ij}^n\big) = T' A$,
 where 
\begin{equation}\label{retrei}
\textstyle
T'=\prod\limits_{j=n-1}^1\,\prod\limits_{i=n}^{j+1} 
\,e_{ij}\big(-{a^j}_{i\rho(j)}({a^j}_{j\rho(j)})^{-1}\big). 
\end{equation}
Note that 
$$
\deg \big (-{a^j}_{i\rho(j)}(a^j_{j\rho(j)})^{-1}\big) \, = \,
\la+\de_{\rho(j)}-\de_i - (\la+\de_{\rho(j)} -\de_j)\, = \,
\de_j - \de_i.
$$
Therefore, in the 
product
for $T'$ each 
$e_{ij}\big (-{a^j}_{i\rho(j)}(a^j_{j\rho(j)})^{-1}\big)
\in \elh$; it is also 
unipotent lower triangular, as $i>j$.
 Hence, $T'$ is homgeneous of 
degree $0$ and is unipotent lower 
triangular.  
Set 
$$
\textstyle
T\,=\,{T'}^{-1}\, =\, 
\prod\limits^{n-1}_{j=1}\,\prod\limits^{n}_{i =j+1} 
\,e_{ij}\big({a^j}_{i\rho(j)}({a^j}_{j\rho(j)})^{-1}\big),
$$ 
which is again a homogeneous unipotent lower triangular matrix 
of degree zero. 
Our construction 
shows that in the matrix $A^{(n-1)}=T^{-1}A$  the leftmost  
non-zero entry in the $i$-th row is 
$a^n_{i\rho(i)}$ which is homogeneous in $\sE$, hence a unit. 
Furthermore, every entry below $a^n_{i\rho(i)}$ is zero. 
The function $\rho$ of the indices is a permutation of 
$\{1, \ldots, n\}$. 
Set 
\begin{equation}\label{hgas}
 U\, = \, \text{diag}(a_{1\,\rho(1)}^n, \ldots, 
a_{n\,\rho(n)}^n), 
%%  U=\begin{pmatrix}
%% c_{1\rho(1)}& 0  & 0 & \dots \\
%%0 & c_{2\rho(2)} & 0 & \dots\\
%% 0 & \dots&  \ddots &\dots \\
%% 0 & \dots& 0  & c_{n\rho(n)}
%%\end{pmatrix}.
\end{equation}
where $\text{diag}(u_1, \ldots, u_n)$ denotes
the  $n\times n$ diagonal matrix with successive
diagonal 
entries $u_1, \ldots, u_n$.
While $U$ need not be homogeneous, its diagonal entries
are all nonzero and homogeneous, hence units of $\sE$;
so, $U$~is invertible in $\sA$.

Clearly $U^{-1}A^{(n-1)}=U^{-1}T^{-1}A$ is a matrix whose leftmost
 non-zero entry in the $i$-th row is the $1$  in the 
$(i,\rho(i))$-position. Furthermore, every entry
below the $(i,\rho(i))$-entry  
is zero. Let $\pi=\rho^{-1}$,
and let $P_{\pi}$ be the permutation matrix of $\pi$. 
Since left multiplication by $P_\rho \ (= P_{\pi}^{-1})$
moves the $i$-th row to the $\rho(i)$-th row
 the matrix  
$$
V\,=\,P_{\pi}^{-1}U^{-1}T^{-1}A
$$ is 
unipotent upper triangular.
 We have 
$A=TUP_{\pi}V$ which we show has the form  
asserted in the proposition.

 As to the homogeneity of these matrices,
we   have seen that $T$ is homogeneous with $\deg(T) = 0$. 
Observe next 
that $U$ and $P_{\pi}$ need not be homogeneous but 
$UP_{\pi}$ is homogeneous.  
% with   $\deg(UP_{\pi})= \deg(A^{(n-1)})=\la$.
 For, $UP_{\pi}$ has its only nonzero entries 
$a^n_{i\rho(i)}$ in the  
 $(i,\rho(i))$-position for  $1\leq i \leq n$
Thus, $UP_\pi$ is obtainable from the homogeneous matrix
$A^{(n-1)}$
by replacing some entries
in 
$A^{(n-1)}$  by $0$'s. 
Hence,  $UP_\pi$ is  homogeneous with $\deg(UP_\pi)
= \deg\big(A^{(n-1)}\big) = \la = \deg(A)$.  Therefore,
$V = (UP_\pi)^{-1}T^{-1}A$ is also homogeneous, 
with $\deg(V) = \deg\big((UP_\pi)^{-1}\big) 
+\deg(T^{-1}) + \deg(A) = 0$.

Next we show that, $P_{\pi} V P_{\pi}^{-1}$ is also 
unipotent upper triangular, so $A=TUP_{\pi}V$ is in
{\it strict}  Bruhat normal form. We have 
\begin{equation}\label{gdgsga}
P_{\pi} V P_{\pi}^{-1}=U^{-1}T^{-1}AP_{\pi}^{-1}.
\end{equation}
Recall the arrangement of entries in the columns of 
$U^{-1}T^{-1}A = U^{-1}A^{(n-1)}$. 
Since 
right multiplication of 
this matrix by $P_\pi^{-1} = P_\rho$
 moves the 
$\rho(i)$-th column 
 to the $i$-th column, $U^{-1}T^{-1}AP_\pi^{-1}$ is 
unipotent upper triangular. Thus,
$P_{\pi} V P_{\pi}^{-1}$ is 
unipotent upper triangular by (\ref{gdgsga}).

It remains  only to show that this strict Bruhat decomposition 
is unique.  
(This uniqueness argument is valid for matrices over any 
ring.)
 Suppose $T_1U_1P_{\pi_1}V_1=T_2U_2P_{\pi_2}V_2$, are
 two strict Bruhat normal forms for the same matrix. Then
\begin{equation}\label{l13}
U_1^{-1}T_1^{-1}T_2U_2=P_{\pi_1}V_1V_2^{-1}P_{\pi_2}^{-1}.
\end{equation}
Since $V_1V_2^{-1}$ is  unipotent upper triangular, we 
can write $V_1V_2^{-1}=\I_n+N$, where 
$\I_n$ is the identity matrix and $N$ is  nilpotent 
upper triangular (i.e., an upper 
triangular matrix with zeros on the diagonal). Note that 
there is no position $(i,j)$ where the  matrices $\I_n$ 
and $N$  both have a non-zero entry. Writing 
\begin{equation}\label{l15}
P_{\pi_1}V_1V_2^{-1}P_{\pi_2}^{-1}\, = \,
P_{\pi_1}P_{\pi_2}^{-1}+P_{\pi_1}NP_{\pi_2}^{-1},
\end{equation} 
the summands on the right again have no overlapping 
nonzero entries.
%the same remains true for the right 
% side of~(\ref{l15}). 
Therefore, as 
$P_{\pi_1}V_1V_2^{-1}P_{\pi_2}^{-1}$ is lower 
triangular by \eqref{l13}, 
each of $P_{\pi_1}P_{\pi_2}^{-1}$ and 
$P_{\pi_1}NP_{\pi_2}^{-1}$ must be lower triangular.
Since $P_{\pi_1}P_{\pi_2}^{-1}=P_{\pi_1\pi_2^{-1}}$ 
is a lower triangular permutation matrix, 
it must be $\I_n$; thus, $\pi_1 = \pi_2$.
Because of the nonoverlapping nonzero entries noted in 
 \eqref{l15}, 
the diagonal entries of $P_{\pi_1}V_1V_2^{-1}P_{\pi_2}^{-1}$
must be $1$'s.  But because the $T_i$ are unipotent lower
triangular and the $U_i$ are diagonal, \eqref{l13} shows that 
the diagonal entries of $P_{\pi_1}V_1V_2^{-1}P_{\pi_2}^{-1}$ coincide
with those of the diagonal matrix
$U_1^{-1}U_2$. Hence, $U_1^{-1}U_2 = \I_n$, 
i.e., $U_1 = U_2$.

Since $\pi_2 = \pi_1$, we can rewrite ~(\ref{l13}) as 
\begin{equation}\label{116}
U_1^{-1}T_1^{-1}T_2U_2\,=\,
P_{\pi_1}V_1P_{\pi_1}^{-1} 
\big(P_{\pi_2}V_2P_{\pi_2}^{-1}\big)^{-1}
\end{equation}
Since the decompositions are strict Bruhat, the right 
side of \eqref{116} is unipotent upper 
triangular while  the left is lower triangular. This 
forces each side to be $\I_n$. Hence,
$V_1=V_2$, $U_1=U_2$ (as we have seeen already),  and 
$T_1=T_2$. This proves the uniqueness.
\end{proof}

\begin{remark} The first part of  the uniqueness 
 proof above (preceding \eqref{116}) shows 
that if $A$ admits a Bruhat decomposition 
$A=TUP_{\pi}V$ (without the  assumption on 
$P_{\pi} V P_{\pi}^{-1}$),  then $\pi$ and $U$ are 
 uniquely determined. 
\end{remark}

\begin{theorem}\label{homde}
Let $\sE$ be a graded division ring. Let 
$\sA=\M_n(\sE)(\overline \de)$ where $\overline \de=(\de_1,\dots, \de_n)$,  
$\de_i \in \Gamma$. 
Then there is a Dieudonnn\'e determinant group 
homomorphism 
$$  \textstyle
\det_\sE\colon \GL_n^h(\sE)(\overline \de) \,\longrightarrow 
\,
\sE_h^*/[\sE_h^*,\sE_h^*].
$$
If $A\in \GL_n^h(\sE)(\overline \de) = \sA_h^*$
has strict Bruhat decomposition $A = TUP_\pi V$
 with ${U = \diag(u_1, \ldots, 
u_n)}$ as in Prop.~\eqref{bruhatf}, then 
\begin{equation}\label{detforBru}
\textstyle 
\det_\sE(A) \, = \, \sgn(\pi) \,u_1\ldots u_n \,[\sE_h^*, \sE_h^*].
\end{equation}
Moreover, if $\det_0\colon\sA_0^* \to \sE_0^*/
[\sE_0^*, \sE_0^*]$ is the Dieudonn\'e determinant for 
the semisimple ring $\sA_0$, then there is a 
commutative diagram
\begin{equation}\label{ddetE0}
\begin{split}
\xymatrix{
\sA^*_0 \ar[r]^{\det_0\  \  \  \  \,}\ar[d] & 
\sE_0^*/[\sE_0^*,\sE^*_0] \ar[d]\\
\sA^*_h \ar[r]^{\det_\sE\  \  \  \  \,\,} & \sE_h^*/[\sE_h^*,\sE_h^*]
}
\end{split}
 \end{equation}
\end{theorem}

\begin{proof}
Throughout the proof we assume that $(\de_1, \ldots,
\de_n) = (\varepsilon_1, \ldots,\varepsilon_1,
\varepsilon_2, \ldots,
\varepsilon_2, \ldots, \varepsilon_k, \ldots,
\varepsilon_k)$  with 
each~$\varepsilon_\ell$ occurring $r_\ell$ times and the 
cosets $\varepsilon_1+\Gamma_\sE, \ldots, 
\varepsilon_k+\Gamma_\sE$
distinct in $\Gamma/\Gamma_\sE$.  There is no loss of 
generality 
with this assumption, in view of \eqref{wdeild}.  Thus, 
any
matrix $B$  in 
$\sA_0$ is in block diagonal form,
say with diagonal blocks $B_1, \ldots , B_k$, with 
each $B_\ell \in M_{r_\ell}(\sE_0)$; we will  
identify
$$
\sA_0 \, = \, \M_{r_1}(\sE_0) \times \ldots \times 
\M_{r_\ell}(\sE_0),
$$
by identifying $B$ with $(B_1, \ldots, B_k)$, which we 
call the block decomposition of $B$.

We first assume that $\sE_0 \ne \mathbb F_2$, the field
with two elements;
  the exceptional case will be treated toward the 
end of the proof.  

It is tempting to use formula \eqref{detforBru}
as the definition of $\det(A)$.  But since it is 
difficult to show that the resulting 
function is a group homomorphism, we take a different 
tack.

We call a matrix $M$ in $\sA$ a {\it monomial matrix}
if $M$ has exactly one nonzero entry in each row and in 
each column, and if each nonzero entry lies in $\sE^*$.    
Clearly, $M$ is a monomial matrix if and only if 
$M = UP$ where $U$ is a diagonal matrix with every
diagonal entry a unit, and $P$ is a permutation matrix.
Moreover, $P$ and $U$ are uniquely determined by $M$.
The set $\mathcal M$ of all monomial matrices in 
$\sA$ is a subgroup of $\sA^*$, and the set 
$\mathcal M^h$
 of all homogeneous monomial matrices is a subgroup 
of $\sA_h^*$.  Define a function 
$$
\Delta\colon \mathcal M^h \,\longrightarrow \sE_h^*/
[\sE_h^*, \sE_h^*] 
$$
by
$$
\Delta(UP_\pi) \, = \ \sgn(\pi)\,u_1 u_2\ldots u_n\, 
[\sE_h^*,\sE_h^*],
\quad \text{where\ \,}U = \text{diag}(u_1, \ldots,
u_n).
$$
This $\Delta$ is clearly well-defined, 
since $UP_\pi$ determines $U$ and 
the permutation matrix $P_\pi$
for $\pi$ in the symmetric group $S_n$. 
Note also that $\Delta$ is a group homomorphism.
For,  if $M = UP_\pi$ and $M'= U' P_{\pi'}$ with 
$U = \text{diag}(u_1, \ldots, u_n)$ and ${U' = 
\text{diag}(u_1', \ldots, u_n')}$, then  
$$
MM' \, = \, \big(UP_\pi U' P_\pi ^{-1}\big) P_{\pi \pi'}
\quad\text{with} \quad
UP_\pi U' P_\pi ^{-1} = \text{diag}
\big(u_1u_{\pi^{-1}(1)}', \ldots, 
u_nu_{\pi^{-1}(n)}' \big).
$$
It follows immediately that $\Delta(MM')=\Delta(M)
\Delta(M')$.  

Recall that 
$\sA_0 = \M_{r_1}(\sE_0) \times \ldots 
\times \M_{r_k}(\sE_0)$. 
Each component $\GL_{r_\ell}(\sE_0)$ of $\sA_0^*$ has
a Dieudonn\'e determinant function  $\det_\ell$ mapping 
it to 
$\sE_0^*/[\sE_0^*, \sE_0^*]$, and these maps are used
to define the Dieudonn\'e determinant 
${\det_0\colon \sA_0^* \to \sE_0^*/[\sE_0^*, \sE_0^*]}$
for the semisimple
ring $\sA_0$ by
\begin{equation}\label{det0}
\textstyle
\det_0(B_1, \ldots, B_k) \,=\, 
\prod\limits_{\ell=1}^k \det_\ell(B_\ell).
\end{equation}
Set $\ov{\det_0}$ to be the composition 
$\sA_0^* \xrightarrow{\det_0} \sE_0^*/[\sE_0^*, \sE_0^*]
\longrightarrow \sE_h^*/[\sE_h^*, \sE_h^*]$.  We claim that 
if $M \in \mathcal M^h$ has degree~$0$, then 
\begin{equation}\label{det0M}
\textstyle
\Delta(M) \,= \,\overline{\det_0}(M). 
\end{equation}
%%%
For, as $M\in \sA_0$, it follows that $U$ and $P$  lie 
in $\sA_0$, and when we view $M=(M_1, \ldots,M_k)$,
${U= (U_1,\ldots,U_k)}$, $P= (P_1,\ldots,P_k)$,
we have:   in each $\ M_{r_\ell}(\sE_0)$,
$U_\ell$ is a diagonal, say $U_\ell = \diag(u_{\ell1}, 
\ldots, u_{\ell r_\ell})$,
 and $P_\ell$ a permutation matrix, say 
$P_\ell = P_{\pi_\ell}$ for some 
$\pi_\ell \in S_{r_\ell}$,  and  $M_\ell = U_\ell P_\ell$.  
So, $M_\ell$ is a monomial 
matrix in $\M_{r_\ell}(\sE_0)$. Since each $M_\ell$ has 
(nonstrict) Bruhat
decomposition $M_\ell = 
\I_{r_\ell} U_\ell P_{\pi_\ell}\I_{r_\ell}$ in 
$\GL_{r_\ell}(\sE_0)$,\cite[\S20, Def.~1, Cor.~1]{draxl} yields
$\det_\ell(M_\ell) = \sgn(\pi_\ell)\,u_{\ell 1}\ldots 
u_{\ell r_\ell} \, [\sE_0^*,\sE_0^*]$.   
Moreover, as  $P = P_\pi$, where $\pi = (\pi_1,\ldots,\pi_k)$ 
when we view
$S_{r_1}\times \ldots \times S_{r_k} \subseteq S_n$, 
we have
$\sgn(\pi) = \sgn(\pi_1) \ldots \sgn(\pi_k)$.
Thus,
$$
\textstyle
\det_0(M)\, = \,\prod\limits
 _{\ell = 1}^k \big(\sgn(\pi_\ell)\,u_{\ell 1}
\ldots u_{\ell r_\ell} \big)\, 
[\sE_0^*, \sE_0^*]
\,=\, \sgn(\pi) \prod\limits_{\ell = 1}^k  (u_{\ell 1}
\ldots u_{\ell r_\ell} )  [\sE_0^*, \sE_0^*],
$$
which yields \eqref{det0M}.

We next claim that every matrix $A$ in $\sA_h^*$ is 
expressible
(not uniquely) in the form $A = CM$, where 
$C \in [\sA_0^*, \sA_0^*]$ and $M \in \mathcal M^h$.  
For this, consider first an elementary matrix 
$e \in \elh$.  The block form of $e$ is 
$(e_1, \ldots, e_k)$, where clearly one $e_\ell$ is an
elementary matrix in $\M_{r_\ell}(\sE_0)$ and all the 
other blocks are identity 
matrices.
Since every elementary matrix in $\M_{r_\ell}(\sE_0)$ lies
in $[\GL_{r_\ell}(\sE_0), \GL_{r_\ell}(\sE_0)]$ by 
\cite[\S20, Th.~3, Th.~4(i)]{draxl} 
(as $\sE_0 \ne \mathbb F_2$ by assumption)
it follows thet $e \in [\sA^*_0,\sA^*_0]$. Now, take 
any $A\in \sA_h^*$, with its strict Bruhat decomposition
$A = TUP_\pi V$ as in Prop.~\ref{bruhatf}.  Then, 
$T$ is a product of elementary matrices in $\sA_0^*$;
so $T\in [\sA^*_0,\sA^*_0]$.  Moreover the transpose $V^t$ 
of $V$ is unipotent lower triangular of degree $0$.  The
unique strict Bruhat normal form of $V^t$  is clearly 
$V^t = V^t \I_nP_{\id}\I_n$.  Hence,  Prop.~
\ref{bruhatf} shows that $V^t$ is a product of 
matrices in $\elh$.  Therefore, $V^t \in
[\sA^*_0,\sA^*_0]$, which implies that $V \in 
[\sA^*_0,\sA^*_0]$.  Now, let $M = UP_\pi \in 
\mathcal M^h$,
and let $C = TMVM\inv = AM\inv$.  Because $V\in 
[\sA^*_0,\sA^*_0]$ and $M$ is homogeneous, 
$M V M \inv\in [\sA^*_0,\sA^*_0]$. (For take any 
$Z_1, Z_2 \in \sA_0^*$. Then, $M[Z_1,Z_2]M\inv
= [MZ_1M\inv ,MZ_2M\inv]\in [\sA^*_0,\sA^*_0]$,
as each $MZ_i M\inv \in \sA_0^*$.) Hence, $C \in 
[\sA^*_0,\sA^*_0]$, so $A = CM$, as claimed.

Define
$\det_\sE\colon \sA_h^* \to \sE_h^* /[\sE_h^*,\sE_h^*]$ by 
$$
\textstyle
\det_\sE(CM) \, = \, \Delta(M), \qquad \text{for any}\  
C \in [\sA_0^*, \sA_0^*], \ M \in \mathcal M^h.  
$$

To see that $\det_\sE$ is well-defined,  
suppose $C_1M_1 = C_2M_2$ with $C_1,C_2\in 
[\sA_0^*, \sA_0^*]$ and $M_1, M_2\in \mathcal M^h$.
Then,
\begin{equation*}\label{MisC}
M_1M_2^{-1} \, = \, C_1^{-1} C_2\, \in [\sA_0^*, \sA_0^*].
\end{equation*}
Hence, $\deg(M_1M_2^{-1}) = 0$ and $\det_0(M_1M_2\inv)
= \det_0(C_1\inv C_2) = 1$, which implies 
that also  ${\overline{\det_0}(M_1M_2\inv)
= 1}$.   So, by \eqref{det0M}
$\Delta(M_1M_2\inv) 
= 1$.  Since $\Delta$ is a group homomorphism, 
it follows that $\Delta(M_1) = \Delta(M_2)$.  Thus,
$\det_\sE$ is well-defined.  To see that it is a group
homomorphism, take any $C,C' \in [\sA^*_0,\sA^*_0]$,
and $M, M' \in \mathcal M^h$. Then, 
$$
(CM)(C'M')\, = \, \big(C(MC'M\inv)\big)(MM').
$$
Since $C' \in [\sA^*_0,\sA^*_0]$, we have
$MC' M\inv \in [\sA^*_0,\sA^*_0]$, as noted above;
so, $C(MC' M\inv) \in [\sA^*_0,\sA^*_0]$.  Also,
$MM' \in \mathcal M^h$.  Hence,
$$
\textstyle \det_\sE\big((CM)(C'M')\big)
\, = \, \Delta(MM') \, = \, \Delta(M)\Delta(M')
\, = \, \det_\sE(CM) \det_\sE(C'M'); 
$$
so, $\det_\sE$ is a group homomorphism.  For $A \in 
\sA^*_h$ with strict Bruhat decomposition
 $A = TUP_\pi V$, we have 
seen that $A = CM$ with $M = UP_\pi\in \mathcal M^h$
and $C = T MV M\inv \in [\sA_0^*,\sA_0^*]$, so 
$\det_\sE(A) = \Delta(M)$, which yields 
formula \eqref{detforBru}.

We now dispose of the  exceptional case 
where $\sE_0 = \mathbb F_2$.  
When this holds, replace $[\sA_0^*,\sA_0^*]$
in the proof by $\sA_0^*$, and the argument 
goes through. 
Observe that now  if $M\in \mathcal M^h$
with $\deg(M) = 0$, then $\Delta(M) = 1$.  For, 
all nonzero entries of $M$ then lie in $\sE_0^* = \{1\}$
and the $\sgn(\pi)$ term in the formula for $\Delta(M)$
drops out as $\chr(\sE_0) = 2$. This replaces
use of \eqref{det0M} in the proof.  There is no need
to invoke $\det_0$, which is in fact trivial here
as $|\sE_0^*| = 1$.  The argument that a homogeneous 
elementary matrix $e$  lies in $[\sA_0^*,\sA_0^*]$
is replaced by the tautology that $e \in \sA_0^*$.

Turning to diagram \eqref{ddetE0}, take any $A\in \sA_0^*$,
with strict Bruhat decomposition 
$A = TUP_\pi V$.  Then, $\det(UP_\pi) = \deg(A) = 0$, 
so $U$ and $P_\pi$ lie in $\sA_0^*$.   Take the block  
decomposition $A = (A_1, \ldots, A_k)$ and likewise for 
$T, U, P, V$.  Then, $P_\pi = (P_{\pi_1}, \ldots,
P_{\pi_k})$, where $\pi = (\pi_1, \ldots, \pi_k)$
when we view $S_{r_1} \times \ldots\times S_{r_k} \subseteq
S_n$.  Note that $A_\ell = T_\ell U_\ell P_{\pi_\ell}
V_\ell$ is the strict Bruhat decomposition of 
$A_\ell$ in $\GL_{r_\ell}(\sE_0)$ for 
$\ell = 1, 2, \ldots, k$. 
So, 
$$
\textstyle
\det_0(A) \, = \, 
\prod\limits _{\ell = 1}^k\det_\ell(A_\ell)
\, = \, \prod \limits _{\ell = 1}^k 
\det_\ell(U_\ell P_{\pi_\ell})
\, = \, \det_0(UP_\pi).
$$
Hence, invoking \eqref{det0M} for $UP_\pi\in \mathcal M^h$,
$$
\textstyle
 \ov {\det_0}(A) \, = \, \ov{\det_0}(UP_\pi)
\, = \, \det_\sE(UP_\pi) \, = \, \det_\sE(A),
$$
showing that diagram \eqref{ddetE0} is commutative.
\end{proof}

In a matrix ring $\M_r(R)$ over any ring $R$, 
for any $a\in R$ we write $\D_r(a)$ for the 
diagonal matrix $\diag(1, \ldots, 1, a)$.

\begin{proposition}\label{kerdetE}
Let $\sA = \M_n(\sE)(\ov\de)$ 
with $(\delta_1, \ldots, \de_n) =
(\varepsilon_1, \ldots,\varepsilon_1,
\ldots, \varepsilon_k, \ldots, \varepsilon_k)$ as in 
the proof of Th.~\ref{homde}.
If $\Gamma_\sE$ is $n$-torsion free, then  
$$
\textstyle
\ker(\det_\sE) \,= \,\big\langle \elh \big\rangle
\cdot \{\big(\D_{r_1}(c_1), \ldots,\D_{r_k}(c_k)\big)
\mid \text{each } c_i \in \sE_0^* \text{ and }
c_1 \ldots c_k \in [\sE_h^*,\sE_h^*]\}.
$$
\end{proposition}

Here, $\elh$ denotes the group of homogeneous elementary matrices, 
as in \eqref{elh}, and $\big(\D_{r_1}(c_1), \ldots,\D_{r_k}(c_k)\big)$
denotes the block diagonal matrix with diagonal blocks
$\D_{r_1}(c_1), \ldots,\D_{r_k}(c_k)$.

\begin{proof}
Suppose $A \in \sA_h^*$ and $\deg(A) = \lambda \ne 0$,
and let $A = TUP_\pi V$ be the strict Bruhat decomposition 
of $A$, with $U = \diag(u_1, \ldots, u_n)$. Since the 
monomial matrix $U P_\pi$  is homogeneous of degree
$\lambda$ with $(i, \pi\inv(i))$-entry $u_i$, we have
$\deg(u_i) = \lambda+ \delta_i - \delta_{\pi\inv(i)}$. 
So $\deg(\sgn(\pi)u_1\ldots u_n) = n\lambda \ne 0$, 
as $\Gamma_\sE$ is $n$-torsion free. But, 
$[\sA_h^*, \sA_h^*]\subseteq \sA_0^*$, as every commutator
of homogeneous matrices has degree $0$.  Hence, 
$\det_\sE(A) \ne 1$. Thus, $\ker(\det_\sE) \subseteq 
\sA_0^*$. 

Note that every homogeneous elementary matrix $e$
has stict Bruhat decomposition $e = e\I_nP_{\id} \I_n$
or $e = \I_n\I_n P_{\id} e$.  In either case, 
$\det_\sE(e) = 1$.  This shows that $\langle
\elh \rangle\subseteq \ker(\det_\sE)$.  

Now take $A \in \sA_0^*$ with block decomposition 
$(A_1, \ldots, A_k)$.  By \cite[\S20, Th.~2]{draxl}, 
each $A_\ell$ is expressible in $\GL_{r_\ell}(\sE_0)$
as $A_\ell = B_\ell \D_{r_l}(c_\ell)$ for some 
$c_\ell \in \sE_0^*$, where $B_\ell$ is a product of 
elementary matrices in $M_{r_\ell}(\sE_0)$.
So, $(\I_{r_1}, \ldots, \I_{r_{\ell-1}}, B_\ell,
  \I_{r_{\ell+1}}, \ldots, \I_{r_k})$ is a product
of the corresponding homogeneous elementary matrices in~ 
$\sA_0$.
Hence, $A = B D$, where $B = (B_1, \ldots, B_k)
 \in \langle \elh\rangle$
and $D = (\D_{r_1}(c_1), \ldots \D_{r_k}(c_k))$, which
is a diagonal matrix in $\sA_0$.  Thus,
$$
\textstyle
\det_\sE(A)  \, = \, \det_\sE(B)\det_\sE(D)
\, = \, c_1\ldots c_k\, [\sE_h^*,\sE_h^*].
$$
So, $A\in \ker(\det_\sE)$ if and only if 
$c_1\ldots c_k \in [\sE_h^*,\sE_h^*]$, which 
yields the proposition.
\end{proof}

Recall that a graded division ring $\sE$ with 
center $\sT$ is said to be {\it unramified} if 
$\Gamma_{\sE}=\Gamma_{\sT}$. 
 In Th.~\eqref{extseq}(iv) below we 
will show that homogeneous $\SK$ of unramified graded 
division algebras is Morita invariant. For 
non-stable $K_1$, we have the following: 

%Here as a corollary of Theorem~\ref{homde}, we show that (nonstable) homogeneous $K_1$ of these algebras is Morita invariant as well. 

\begin{cor}\label{homdeun} 
Let $\sE$ be a graded division ring and
let $\sA=\M_n(\sE)$ 
with unshifted grading. Suppose $\Ga_\sE$ is $n$-torsion
free, $\sE$ is unramified, and $\M_n(\sE_0) 
\not = \M_2(\mathbb F_2)$. Then 
 $\det_\sE$  induces a group monomorphism   
\[ \GL_n^h(\sE)\big/ [\GL_n^h(\sE), \GL_n^h(\sE)] \ 
\hookrightarrow \ \sE_h^*\big/[\sE_h^*,\sE_h^*].\]
\end{cor}
\begin{proof}
We need to show that $\ker(\det_\sE) = 
[\GL_n^h(\sE), \GL_n^h(\sE)]$.  The inclusion $\supseteq$
is clear as $\det_\sE$ is a group homomorphism 
mapping into an abelian group.  For the reverse inclusion,
note that as $\sA_0$ is simple, Prop.~\ref{kerdetE} says 
$$
\textstyle
\ker(\det_\sE) \, = \, \big\langle \elh\big\rangle\cdot 
\big\{\D_n(a) \mid a\in [\sE_h^*, \sE_h^*]\big\}.
$$
Because $\Ga_\sE = \Ga_\sT$ where $\sT$ is the center of 
$\sE$, we have $\sE_h^* = \sT_h^* \cdot \sE_0^*$, hence 
$[\sE_h^*, \sE_h^*] =[\sE_0^*, \sE_0^*]$.  Thus, 
$$
\big\{\D_n(a) \mid a\in [\sE_h^*, \sE_h^*]\big\}
 \, =\, [\D_n(\sE_0^*), \D_n(\sE_0^*)]
\, \subseteq \,[\GL_n^h(\sE), \GL_n^h(\sE)].
$$
Also, as $\sA_0 = \M_n(\sE_0)$, the homogeneous 
elementary matrices of $\sA$, which all have degree $0$, 
are the same as the elementary matrices of 
$\M_n(\sE_0)$;  since $\sE_0$ is a division ring, 
by \cite[\S20, Th.~4, Lemma~4]{draxl}
these all lie in $[\GL_n(\sE_0), \GL_n(\sE_0)]
\subseteq [\GL_n^h(\sE), \GL_n^h(\sE)]$, as 
$\M_n(\sE_0) \ne \M_2(\mathbb F_2)$.  
Hence, $\ker(\det_\sE) \subseteq 
[\GL_n^h(\sE), \GL_n^h(\sE)]$, completing the proof.
\end{proof}

\section{Homogeneous $\SK$}

Throughout this section we consider graded division 
algebras $\sE$, i.e., $\sE$ is a graded division ring 
 which is finite-dimensional as a graded vector space 
over its center $\sT$. In addition, we assume that the abelian
group $\Gamma$ containing $\Gamma_\sE$ 
is torsion free. The assumption on $\Gamma$ implies that
every unit in $\sE$ is actually homogeneous, so $\sE_h^* 
= \sE^*$. This assumption 
also implies that $\sE$ has no zero divisors. 
(These properties  follow easily from the fact that the 
torsion-free
 abelian group $\Gamma_\sE$ can be made into
a totally ordered group, see, e.g. \cite[p.~78]{hwcor}.) Hence, 
$\sE$ has a quotient division ring 
obtained by central localization, $q(\sE) = 
\sE \otimes_\sT q(\sT)$, where $q(\sT)$~is the quotient
field of the integral domain $\sT$.
In addition, every graded 
module $\sM$
over  $\sE$ is a free module with well-defined rank; we 
thus call $\sM$ a graded vector space over $\sE$, 
and write $\dim_\sE(\sM)$ for $\text{rank}_\sE(\sM)$.  
This applies also for graded modules over $\sT$, which is 
a commutative graded division ring.  We write $[\sE:\sT]$
for $\dim_\sT(\sE)$, and $\ind(\sE) = \sqrt{[\sE:\sT]}$.
Clearly, $[\sE: \sT] = [q(\sE): q(\sT)]$,
so $\ind(\sE) = \ind(q(\sE)) \in \mathbb N$.    
Moreover, in (\cite[Prop.~5.1]{boulag}
 and~\cite[Cor.~1.2]{hwcor} it was observed that $\sE$ is 
an Azumaya algebra over $\sT$. 

In general for an Azumaya algebra $A$ of constant rank 
$m^2$
over a commutative ring $R$, there is a  commutative 
ring $S$ faithfully flat over $R$ which splits $A$, i.e., 
$A\otimes_R S \cong \M_m(S)$. For $a \in A$, 
considering $a \otimes 1$ as an element of $\M_m(S)$, 
one then defines the {\it reduced characteristic 
polynomial}, $\chr_A(X,a)$, 
 the {\it reduced trace}, $\Trd_A(a)$, and the 
{\it reduced norm}, $\Nrd_A(a)$, of $a$ by
  $$
\chr_A(X,a) \ = \ \det(X\I_m-(a\otimes1)) \ = \ 
X^m-\Trd_A(a)X^{m-1}+\ldots+(-1)^m\Nrd_A(a).
$$
in the polynomial ring $S[X]$.
Using descent theory, one shows that $\chr_A(X,a)$ is 
independent of $S$ 
and of the choice of $R$-isomorphism 
$A\otimes_RS \cong \M_m(S)$, and 
that $\chr_A(X,a)$ lies in $R[X]$;  
furthermore,  the element $a$ is invertible in $A$ if 
and only if 
$\Nrd_A(a)$ is invertible in $R$ (see 
Knus \cite[III.1.2]{knus} 
and Saltman \cite[Th.~4.3]{saltman}). Let $A^{(1)}$ 
denote the multiplicative group 
of elements of $A$ of reduced 
norm $1$. One then  defines the {\it reduced Whitehead 
group} 
of $A$ to be $\SK(A)=A^{(1)}/A'$, where $A'=[A^*,A^*]$ 
denotes the commutator 
subgroup of the group $A^*$ of  units of 
$A$. For any integer $n\ge 1$, the matrix ring
$\M_n(A)$ is also an Azumaya algebra over 
$R$.  One says that $\SK$ is {\it Morita invariant}
for $A$ if 
$$
\SK(\M_n(A)) \cong \SK(A) \quad\text{ for all } n\in \nn.
$$  

Specializing to the case of a graded 
division algebra $\sE$  and the graded   matrix
algebra ${\sA=\M_n(\sE)(\overline \delta)}$, 
where $\ov \de = (\delta_1, \ldots, \de_n) \in \Gamma^n$, we have
 the { reduced Whitehead group} 
\begin{equation} \label{azuski}
 \SK(\sA)\,=\,\sA^{(1)}/[\sA^*,\sA^*],
\text{  where  }\sA^{(1)} 
=\big \{x\in \sA^* \mid \Nrd_{\sA}(x)=1\big\}. 
 \end{equation}
Here $\sA^*$ is the group of units of the 
ring 
$\M_n(\sE)$ (thus the shifted grading on
$\sA$ does not 
affect $\SK(\sA)$).  Restricting to the homogeneous 
elements 
of $\sA$ we define  
\begin{equation}\label{zuhgoi}
\hSK(\sA)= \sA_h^{(1)}\big/[\sA_h^*,\sA_h^*]\ \ 
\text{where}\ \ \sA_h^{(1)} = 
\big \{x\in \sA_h^* \mid \Nrd_{\sA}(x)=1\big\}.
\end{equation} 
 To distinguish these two groups, we call the second one 
the {\it homogeneous reduced Whitehead group}
of $\sA$. 
These groups coincide for $n=1$, i.e, 
$\hSK(\sE) = \SK(\sE)$.  For, $\sE^* = \sE^*_h$, 
as noted above.
(See \cite{HW1} for an extensive study 
of $\SK$ of graded division algebras.)  

The question naturally arises  whether $\SK(A)$ is Morita 
invariant for an Azumaya algebra $A$.
When $A$~is a  central simple algebra
this is known to be the case (see, e.g.,
\cite[\S22, Cor.~1]{draxl} or  
\cite[\S16.5, Prop.~b]{pierce}). We will answer 
the analogous question for  homogeneous reduced 
Whitehead groups when $A$ is a graded division algebra $\sE$
by establishing an exact sequence 
relating  $\hSK(\M_n(\sE))$  and $\SK(\sE)$ 
(Th.~\ref{extseq}) and producing  examples showing 
that  they sometimes differ (Ex.~\ref{trunshift}); thus, 
$\hSK$ is not Morita invariant.  We will see in fact that, 
as $n$ varies, $\hSK(\M_n(\sE))$ depends only on the congruence
class of $n$ modulo a constant $e$ dividing the ramification
index of $\sE$ over its center.  Furthermore, $\hSK(\M_n(\sE))
\cong \SK(\sE)$ whenever $n$ is prime to $e$.

A major reason why $\SK^h(\sA)$ is more tractable than
$\SK(\sA)$  for $\sA = \M_n(\sE)(\ov \de)$ is that $\sA_h^{(1)}$ consists of homogeneous 
elements of degree $0$, as we next show.   This will
allow us to use 
the Dieudonn\'e determinant for the semisimple algebra~ 
$\sA_0$ to relate $\SK^h(\sA)$ to $\SK(\sE)$.

\begin{lemma}\label{indeg0}
With the hypotheses on $\sE$ as above, let 
$\sA = \M_n(\sE)(\ov \delta)$ for $\ov \de 
= (\de_1, \ldots, \de_n) \in 
\Gamma^n$.  Let $\sT$ be the center of $\sE$. Then,
$\Nrd_\sA(\sA_\la) \subseteq \sT_{ns\la}$ for any
 $\la \in \Gamma_\sA$, where $s = \ind(\sE)$.  Hence,
$\sA_h^{(1)}\subseteq \sA_0^*$.
\end{lemma}
  
\begin{proof}
For caluclating $\Nrd_\sA$, we split $\sE$ using a 
graded faithfully flat extension of its center 
$\sT$, in order to preserve the graded structure.
For this we employ a maximal graded subfield 
$\sL$ of $\sE$.  Associated to the graded field $\sT$
there is a graded Brauer group $\gbr(\sT)$ of equivalence 
classes of graded division algebras with center $\sT$.
See \cite{hwcor,tigwad1} for properties of graded Brauer 
groups. In particular, there is a commutative diagram 
of scalar extension homomorphisms,
 \begin{equation*}
\xymatrix{
 \gbr(\sT) \ar[r]  \ar@{^{(}->}[d]^{- \otimes_\sT q(\sT)}& 
\gbr(\sL) \ar@{^{(}->}[d]^{- \otimes_L q(\sL)}\\
\br(q(\sT)) \ar[r] & \br(q(\sL)),
}
 \end{equation*}
where the vertical maps are injective.
If $\sL$ is a maximal graded subfield of 
$\sE$, then $[\sL:\sT]=\ind(\sE)$
by the graded Double Centralizer Theorem \cite[Prop.~1.5]{hwcor}. 
Since 
$[q(\sL):q(\sT)]=[\sL:\sT]=\ind(\sE)=\ind(q(\sE))$,  
it follows that $q(\sL)$ is a maximal subfield of the 
division ring $q(\sE)$, which is  known to be a 
splitting field for~$q(\sE)$ (see \S9, Cor.~5 
in~\cite{draxl}). The commutativity of the diagram above
and the injectivity of vertical arrows imply that $\sL$ 
 splits  $\sE$ as well, i.e.,  $\sE\otimes_{\sT} \sL 
\cong_{\gr} \M_s(\sL)(\overline {\gamma})$, for some 
 $\overline{\ga}=
(\ga_1,\dots,\ga_s)\in \Gamma^s$, where $s = \ind(\sE)$.
Moreover $\sL$ is a free, hence faithfully flat, $\sT$-module.

The graded field $\sL$ also splits $\sA = \M_n(\sE)(\overline \de)$,
where $\overline \de=(\de_1,\dots,\de_n) \in \Gamma^n$. Indeed, 
 $$
\sA \otimes_{\sT} \sL\cong_{\gr} \M_n(\sE)
(\overline \de)  \otimes_{\sT} \sL \,\cong_{\gr} \, 
\M_n(\sE\otimes_{\sT} \sL)(\overline \de) \,
\cong_{\gr}\, \M_n\big(\M_s(\sL)(\overline \ga)\big)
(\overline \de) \,\cong_{\gr}\, \M_{sn}(\sL)
(\overline{\omega}),
$$ 
where $\overline{\omega}=
(\ga_i+\de_j)$, $1\leq i \leq s$, $1\leq j \leq n$.
For a homogeneous element $a$ of $\sA$ 
with $\deg(a) = \la$, its image $a\otimes 1$ in 
$\sA \otimes _\sT \sL$ is also homogeneous of degree 
$\lambda$, and $\Nrd_\sA(a) = \det(a\otimes 1)$. 
But, as noted in (\ref{jsjgfd}) above, 
$\det(s\otimes1) \in \sT_{ns\lambda}$.  
Thus, $\Nrd(\sA_\lambda) \subseteq \sT_{ns\la}$.
If $\Nrd_\sA(a) = 1 \in \sT_0$, then $\deg(a) = 0$,
as $\Gamma$ is assumed torsion free.
Thus, $\sA^{(1)}_h \subseteq \sA_0$.  
\end{proof}

In order to establish a 
connection between the homogeneous $\hSK(\sA)$  and 
$\SK(\sE)$ 
 we need to relate the reduced norm of $\sA$ to 
that of  $\sA_0$, which we do in the next lemma. 
 Recall that $\sA_0$ is a semisimple ring 
(see (\ref{urnha})). 
For a division algebra $D$, one defines the reduced norm 
map on a semisimple algebra $\M_{r_1}(D)\times \dots
\times \M_{r_k}(D)$ finite-dimensional over its center
as the product of reduced norms of the simple factors.

\begin{lemma} \label{nomng}
With the hypotheses on the graded division algebra 
$\sE$ as above, let $\sA = \M_n(\sE)(\ov \de)$
for ${\ov \de = (\de_1, \ldots , \de_n) \subseteq \Ga^n}$.
Let $\sT$ be the center of $\sE$.
Then,  for $a\in \sA_0$
\begin{equation}\label{nrdnrd}
\Nrd_\sA(a)\,=\,N_{Z(\sE_0)/\sT_0}(\Nrd_{\sA_0}(a))^d,
\end{equation} 
where 
$d=\ind(\sE)/\big (\ind(\sE_0)\,[Z(\sE_0):\sT_0]\big )$. 
\end{lemma}

Here $Z(\sE_0)$ denotes the center of $\sE_0$, which 
is a field finite-dimensional and abelian Galois over 
$\sT_0$.  Also, $N_{Z(\sE_0)/\sT_0}$ denotes the field
norm from $Z(\sE_0)$ to $\sT_0$.

\begin{proof}
After applying a graded isomorphism, we may assume 
that $(\de_1, \ldots , \de_n ) $ has the form\break  $ 
(\varepsilon_1, \ldots, \varepsilon_1, \varepsilon_2,
\ldots , \varepsilon_2, \ldots, \varepsilon_k, 
\ldots, \varepsilon_k)$ as in \eqref{wdeild} above. 
Then,   
$\sA_0 =\M_{r_1}(\sE_0)\times \ldots \times
 \M_{r_k}(\sE_0)$.
Let\break ${a = (a_1, ... , a_k) \in \sA_0}$ with each 
$a_i \in \M_{r_i}(\sE_0)$.  That is, $a$ is in block 
diagonal 
form with diagonal blocks $a_1, \ldots, a_k$;
so, 
$\Nrd_{\sA_0}(a) = \prod_{i=1}^k \Nrd_{\M_{r_i}}(\sE_0)
(a_i)$.    
We thus need to prove that:
\begin{equation}\label{zuzu}
\textstyle
\Nrd_{\sA}(a) \,=\,\prod\limits_{i = 1}^k N_{Z(\sE_0)/\sT_0}
(\Nrd_{\M_{r_i}(\sE_0)}(a_i))^d,
\end{equation}
  Formula~\eqref{nrdnrd} is known for $n= 1$,
i.e., $\sA =\sE$, by \cite[Prop.~3.2]{HW1}.  
The  further 
fact needed here is that
for any  matix  $b$ in $\M_n(\sE)$ in  block 
triangular form, say with
diagonal blocks $b_1,\dots, b_m$, where 
$b_j \in \M_{t_j}(\sE)$, and $t_1+ \ldots +t_m = n$,
we have
\begin{equation}\label{pa22}
\textstyle
\Nrd_{\M_n(\sE)}(b)  \,=\, \prod\limits_{j = 1}^m
\Nrd_{\M_{t_j}(\sE)}(b_j).
\end{equation}
Indeed, if we split $\sE$ by extending scalars, say   
$\sE \otimes_\sT \sL \cong \M_s(\sL)$ for some 
graded field $\sL$,
then \[{\M_n(\sE)\otimes _\sT \sL\cong \M_{ns}(\sL)};\]
  the   matrix for $b\otimes 1$ is again in block 
triangular form with its
diagonal blocks coming from the splitting of the diagonal
 blocks
of $b$.  So formula~\eqref{pa22} follows from the 
determinant formula for
matrices in block triangular form.

Formula~\eqref{pa22}  applied to the block 
diagonal matrix $a$  
shows that it suffices to 
verify that
\begin{equation}\label{2.7}
\textstyle
\Nrd_{\M_{r_i}(\sE)}(a_i) \,=\, N_{Z(\sE_0)/\sT_0}
(\Nrd_{\M_{r_i}(\sE_0)}(a_i))^d
\end{equation}
for each $i$. 
Formula~\eqref{2.7} is 
clearly multiplicative in $a_i$. Moreover,
it holds for any triangular matrix in $\M_{r_i}(\sE_0)$ by
\eqref{pa22} with $t_1 = \ldots = t_m = 1$
and $m=r_i$, since it holds when $\sA = \sE$.
But, we can
always write $a_i  =  e_{i1} c_i e_{i2}$,
where $e_{i1},e_{i2}$ are products of elementary matrices
in $\M_{r_i}(\sE_0)$ and $c_i$ is a diagonal matrix.  This
 is just another way of saying that we can diagonalize 
$a_i$ in $\M_{r_i}(\sE_0)$ 
by elementary row and column operations.  
Thus, formula \eqref{2.7} holds for $a_i$
because $a_i$ is a product of triangular matrices.
This yields~\eqref{nrdnrd}.
\end{proof}

In producing the first  examples of division  
algebras $D$ with nontrivial reduced Whitehead groups, 
Platonov worked in \cite{platonov} with division algebras 
over twice iterated 
Laurent series over a global field.  Ershov later in 
\cite{ershov} generalized and systematized Platonov's approach,
by working with division algebras over arbitrary Henselian 
valued fields.  Ershov encapsulated his results in a 
commutative diagram with exact rows and 
columns which related $\SK(D)$ to various quantities involving
the residue division algebra $\ov D$ and the value group~
$\Gamma_D$ for the valuation on $D$.  More recently it 
was shown in \cite[Th.~3.4]{HW1} that there is a commutative 
diagram
analogous to Ershov's for computing $\SK(\sE)$, where 
$\sE$ is a graded division algebra.  It was also shown
in \cite[Th.~4.8]{HW1}
that Ershov's results for $D$ over a Henselian field 
could be deduced from the corresponding graded ones 
by proving that $\SK(D)\cong \SK(\gr(D))$, where 
$\gr(D)$ is the associated graded division algebra 
of the valued division algebra~$D$. The diagram for 
$\SK(\sE)$ is the vertical $\sE$-plane in the following
diagram \eqref{ershovsk1}.

%%%
\begin{equation}\label{ershovsk1}
\begin{split}
\xymatrix@=10pt{
&    & &&   \scriptstyle{1}  \ar[dd] & \\
& &  &  \scriptstyle{1} \ar[dd]& \\
&    & \scriptstyle{\SK(\sE_0)} \ar@{->}'[r][rr] && 
\scriptstyle{\ker \widetilde N_{\sE}/[\sE_0^*,\sE^*] }
\ar[rr]^-{\Nrd_{\sE_0}} \ar@{->}'[d][dd] && 
\scriptstyle{\widehat H^{-1}(G,\Nrd_{\sE_0}(\sE_0^*))}
\ar[rr] & &  \scriptstyle{1}  \\
& \scriptstyle{\SK(\sA_0)} \ar[rr] 
\ar[ru]^{\scriptscriptstyle{\cong\ \ }} && 
\scriptstyle{\ker \widetilde N_{\sA}/[\sA^*_h,\sA^*_0] }
\ar[rr] \ar[dd] \ar[ru]^{\scriptscriptstyle{\cong\ \ }} && 
\scriptstyle{\widehat H^{-1}(G,\Nrd_{\sA_0}(\sA_0^*))} 
\ar[r] \ar[ru]^{\scriptscriptstyle{\cong\ \ }} &   
\scriptstyle{1}  \\
&   &  \scriptstyle{ \Gamma_{\sE} /\Gamma_{\sT} \wedge 
\Gamma_{\sE} /\Gamma_{\sT} } %[E^*,E^*]/[E_0^*,E^*] 
\ar@{->}'[r][rr] & & \scriptstyle{ \sE^{(1)}/[\sE_0^*,
\sE^*] }
\ar[rr]  \ar@{->}'[d][dd]&& 
\scriptstyle{ \SK(\sE)} \ar[r] &  \scriptstyle{1} \\
   &  \scriptstyle{ \Gamma_{\sA} /\Gamma_{\sT} \wedge 
\Gamma_{\sA} /\Gamma_{\sT} } \ar[ru]^{\eta_n\ } 
\ar[rr] &&  \scriptstyle{ \sA_h^{(1)}/[\sA^*_h,\sA^*_0] }
\ar[rr]  \ar[ru]^{\scriptscriptstyle{\cong\ \ }} 
\ar[dd]^{\widetilde N}&& 
\scriptstyle{ \hSK(\sA)} \ar[ru] \ar[r] &  
\scriptstyle{1} \\
&  & & & \scriptstyle{ \mu_\dlambda(\sT_0) \cap  
\widetilde N (\sE_0^*)} \ar[dd]&\\
 & & &\scriptstyle{ \mu_\dlambda(\sT_0) \cap  
\widetilde N (\sA_0^*)} 
\ar[dd] \ar[ru]^{\scriptscriptstyle{\cong\ }}&\\
&  & & &  \scriptstyle{1}  &\\
 & & &  \scriptstyle{1}  &}
 \end{split}
 \end{equation} 
%%%

 This 
diagram  shows the close connections
between $\SK(\sE)$ and $\hSK(\sA)$, where 
$\sA = \M_n(\sE)$ with unshifted grading.
The diagram is  commutative  with exact rows 
and columns. 
 The group $G$ appearing there is 
${G = \Gal(Z(\sE_0)/\sT_0)}$, where $\sT$ is the center
of $\sE$, and $Z(\sE_0)$ is the center of $\sE_0$;
it is known that $Z(\sE_0)$ is Galois over 
$\sT_0$, and that $G$ is a homomorphic image
of $\Gamma_\sE/\Gamma_\sT$, so $G$ is abelian.  
Also, $d =
 \ind(\sE)\big/\big(\ind(\sE_0)\,
[Z(\sE_0):\sT_0]\big)$, and  $\mu_d(\sT_0)$ is 
the group of those
$d$-th roots of unity lying in $\sT_0$.
The map $\widetilde N_{\sE}$ is the composition 
$\widetilde N_{\sE} = N_{Z(\sE_0)/\sT_0}\circ 
\Nrd_{\sE_0}\colon
\sE_0^* \to \sT_0^*$;  the map $\widetilde N_{\sA}$ 
is defined analogously. 
Exactness of the rows and column in the vertical
$\sE$-plane is proved in  
\cite[Th.~3.4]{HW1}; exactness in the 
$\sA$-plane is proved analogously, as the 
reader can readily verify.   
The maps from the $\sA$-plane to the $\sE$-plane are 
induced by the Dieudonn\'{e} determinant 
$\det_{\sA_0}$ from~ 
$\sA_0 = \M_n(\sE_0)$ to $\sE_0$. By Lemma~\ref{indeg0}, 
$\sA_h^{(1)} \subseteq \sA_0^*$; moreover, the images of 
$\sA_h^{(1)}$, $[\sA_h^*,\sA_h^*]$ and 
$[\sA_h^*,\sA_0^*]$ in $\sA_0^*/[\sA_0^*,\sA_0^*]$,  
under $\det_{\sA_0}$, lie in the images of 
$\sE^{(1)}$, $[\sE^*,\sE^*]$ and $[\sE^*,\sE_0^*]$ in 
$\sE^*_0/[\sE^*_0,\sE^*_0]$, respectively--- 
see  Prop.~\ref{totalfilt} below, which yields the middle
isomorphism in the lower horizontal plane of the 
diagram. 
Here, $\Gamma_\sA = \Gamma_\sE$ since the grading 
on $\sA$ is unshifted, and the map~$\eta_n$ on the left
 is $x\mapsto nx$.  This diagram gives some insight 
into where to look for differences between 
$\hSK(\sA)$ and $\SK(\sE)$;  the differences are
delineated in Th.~\ref{extseq} below.

Let $\sA = \M_n(\sE)$, with unshifted grading. 
We have the filtration of commutator groups
\begin{equation*}
 [\sA_0^*,\sA_0^*] \, \subseteq \,[\sA_h^*, \sA_0^*]
\, \subseteq \, [\sA_h^*, \sA_h^*] \,\subseteq
\, \sA_h^{(1)},  
\end{equation*}
with $\SK^h(\sA) = \sA_h^{(1)}/[\sA_h^*,\sA_h^*]$.
We  relate the factors in this filtration to 
the corresponding ones for $\sE$ in order
to relate $\SK^h(\sA)$ to $\SK(\sE)$:

\begin{proposition}\label{totalfilt}
Let $\sA = \M_n(\sE)$ with unshifted grading, and suppose
$\sA_0 \ne \M_2(\mathbb F_2)$.  Then,
\begin{equation}\label{factoriso}
\sA^{(1)}\big/\, [\sA_0^*,\sA_0^*] \, \cong \, 
\sE^{(1)}\big/ \, [\sE_0^*, \sE_0^*],
\end{equation}
and this isomorphism maps 
$[\sA_h^*, \sA_0^*]\big/[\sA_0^*, \sA_0^*]$
onto $[\sE_h^*, \sE_0^*]\big/[\sE_0^*, \sE_0^*]$.
\end{proposition}

\begin{proof}
Let 
$\overline{\sA^{(1)}_h}=\sA^{(1)}_h/[\sA_0^*,\sA_0^*]$
and $\overline{\sE^{(1)}}=\sE^{(1)}/[\sE_0^*,\sE_0^*]$.
Note that $\sA_0 = \M_n(\sE_0)$, since the grading 
on $\sA$ is unshifted.
There is a homomorphism 
$$
\eta\colon\overline{\sE^{(1)}} \rightarrow 
\overline{\sA^{(1)}_h} \;  \text{ induced by }
c \mapsto \diag(c,1,1,\ldots, 1).
$$ 
This $\eta$ is well-defined, as 
$\Nrd_\sA(\diag(c,1, \ldots,1)) = \Nrd_\sE(c)$.
Moreover, $\eta$ is surjective, as\break 
$\sA_0^* = \diag(\sE_0^*, 1, \ldots ,1)[\sA_0^*,
\sA_0^*]$ (see \cite[\S22, Th.~1]{draxl}) since $\sA_0 \ne 
\M_2(\mathbb F_2)$.  To get a map in the other 
direction we use the Dieudonn\'e determinant for 
$\sA_0$, 
$$
\textstyle
\det_{\sA_0} \colon \sA_0^* \,\longrightarrow\,
\sE_0^*\big/\,[\sE_0^*,\sE_0^*].
$$
Recall (see \cite[\S22, Th.~1]{draxl}) that $\det_{\sA_0}$
is compatible with reduced norms, i.e. 
$\Nrd_{\sA_0}(a) = \ov{\Nrd_{\sE_0}}(\det_{\sA_0}(a))$
for all $a\in \sA_0^*$, where $\ov{\Nrd_{\sE_0}}
\colon \sE_0^*/[\sE_0^*,\sE_0^*] \to Z(\sE_0)^*$ 
is induced by $\Nrd_{\sE_0}$.
Therefore, if $a\in \sA^{(1)}$, then 
$a\in \sA_0^*$ by Lemma~\ref{indeg0}, so
by Lemma~\ref{nomng} (used for $\sA$ then for $\sE$),
$$
1\,=\,\Nrd_\sA(a)\,=\,
N_{Z(\sE_0)/\sT_0}(\Nrd_{\sA_0}(a))^d\,=\,
N_{Z(\sE_0)/\sT_0}\big(\,\ov{\Nrd_{\sE_0}}(\text{det}_{\sA_0}(a))
\big)^d
\,=\,\ov{\Nrd_{\sE}}(\text{det}_{\sA_0}(a)).
$$
This shows that there is a well-defined homomorphism
$$
\textstyle
\xi\colon \ov {\sA_h^{(1)}}\, \longrightarrow \, 
\ov{\sE_h^{(1)}}\quad \text{induced by $\det_{\sA_0}$}.
$$
Since $\det_{\sA_0}(\diag(c, 1,\ldots,1)) = c\,[\sE_0^*,
\sE_0^*]$, for $c\in \sE_0^*$ we have $\xi \eta = \id$.  Therefore, as
$\eta$ is surjective, $\eta$ and $\xi$ are isomorphisms,
proving \eqref{factoriso}.

Let $\overline{[\sA_h^*,\sA_0^*]} =
[\sA_h^*, \sA_0^*]\big/[\sA_0^*, \sA_0^*]$
and $\overline{[\sE^*,\sE_0^*]} = 
[\sE_h^*, \sE_0^*]\big/[\sE_0^*, \sE_0^*]$.  It remains to 
show that these groups are isomorphic via $\xi$.

Since $\Gamma_\sA = \Gamma_\sE$ as the grading on $\sE$ is 
unshifted, we have $\Ga_\sA^* = \Ga_\sE$.  That is, for
any $s\in \sA_h^*$ there is $e\in \sE^*$ with $\deg(e)
= \deg(s)$.  Then, $s = [s(e\inv\I_n)]\, e\I_n$
with $\deg(s(e\inv \I_n)) = 0$.  Thus, 
$\sA_h^* = (\sE^*\I_n)\sA_0^*$. Recall the general 
commutator identity
\begin{equation}\label{commid}
[ab,c]\,=\,[{}^{a}b,{}^{a}c][a,c]
\qquad\text{where \
${}^{a}x=axa^{-1}$}.
\end{equation} 
Since  $\sA_0^*$ is a normal subgroup of $\sA_h^*$, 
this identity  shows that 
$\overline{[\sA_h^*,\sA_0^*]}$ is generated by 
the images of 
commutators of the form $[c\I_n,a]$, where $c\in \sE^*$ 
and $a\in \sA_0^*$. 
Now  if $\varphi$ is any ring automorphism of $\sE_0$,
then $\varphi$ induces an automorphism of 
$\sA_0 = \M_n(\sE_0)$, again called~$\varphi$, and 
also an automorphism $\ov \varphi$ of $\sE_0^*/[\sE_0^*,\sE_0^*]$.
Because $\varphi$ is compatible with strict Bruhat 
decompositions of matrices, it is compatible with 
$\det_{\sA_0}$, i.e., $\det_{\sA_0}(\varphi(s)) 
= \ov \varphi( \det_{\sA_0}(s))$ for any $s\in \sA_0^*$. 
By applying this to   the automorphism of $\sE_0$
given by conjugation by $c \in \sE^*$, we obtain, for 
any $a\in \sA_0^*$,   
$$
\textstyle
\det_{\sA_0}([c\I_n,a])\,=
\,\det_{\sA_0}(c\I_n \,a\, c^{-1} \I_n)
\det_{\sA_0}(a^{-1})\,
=\, cdc\inv d\inv \,[\sE_0^*,\sE_0^*]\quad\text{ where }
\det_{\sA_0}(a) = d \,[\sE_0^*,\sE_0^*].
$$ 
This shows that $\xi\big(\overline{[\sA_h^*,\sA_0^*]}\big)
=\overline{[\sE^*,\sE_0^*]}$, and hence 
$\eta\big(\overline{[\sE^*,\sE^*_0]}\big)=
\overline{[\sA_h^*,\sA_0^*]}$.
\end{proof}

\begin{theorem}\label{extseq}
Let $\sE$ be a graded division algebra finite-dimensional over
its center 
$\sT$ $($with  $\Ga_\sE$ torsion-free$)$.
For  $n \in \mathbb N$
let $ \sA =\M_n(\sE)$, with unshifted grading,
and assume $\M_n(\sE_0) \ne \M_2(\mathbb F_2)$. Then 
there is an exact sequence
\begin{equation}\label{mainthwxse}
 0 \,\longrightarrow \,[\sE^*,\sE^*]\,\big 
/\big([\sE^*,\sE^*]^n\,[\sE^*,\sE_0^*]\big)
\, \longrightarrow \,
\hSK(\M_n(\sE))\, \xrightarrow{\ \ov\xi \ } \,\SK(\sE) 
\,\longrightarrow\, 0,
 \end{equation}
where $\ov \xi$ is induced by the Dieudonn\'e determinant
$\det_{\sA_0}\colon \sA_0^* \longrightarrow 
\sE^*/[\sE_0^*, \sE_0^*]$.~Furthermore, let\break
 ${\Lambda = \Gamma_{\sE} /\Gamma_{\sT} 
\wedge \Gamma_{\sE}
 /\Gamma_{\sT}}$, a finite abelian group. and let 
$e$ be the exponent of  $\Lambda$. Then,
 \begin{enumerate}
\item [(i)]  
The group $[\sE^*,\sE^*]\big /
[\sE^*,\sE_0^*]$ is a homomorphic image of $\Lambda$. 
Hence,
 $[\sE^*,\sE^*]\,\big /\big([\sE^*,\sE^*]^n[\sE^*,\sE_0^*]\big)$ 
is a homomorphic image of $\Lambda /n\Lambda$.
\smallskip
\item[(ii)]   As $n$ varies, 
$[\sE^*,\sE^*]\big /[\sE^*,\sE^*]^n[\sE^*,\sE_0^*]$ 
depends only on the congruence class of 
$n$ mod $e$. 
\smallskip

\item[(iii)] If $\gcd(n,e) = 1$, then 
$\hSK(\M_n(\sE)) \cong \SK(\sE)$.  
This holds for all $n$ if $\Lambda$ is trivial,
which occurs, e.g. if $\Gamma_{\sE}=\mathbb Z$ or
more generally if $\Gamma_{\sE} /\Gamma_{\sT}$ 
is cyclic.

\smallskip
  
\item [(iv)] If $\sE$ is unramified over $\sT$, 
then $\hSK(\M_n(\sE)) \cong \SK(\sE)\cong\SK(\sE_0)$. 

\smallskip

\item [(v)] Suppose $\sE$ is totally ramified over $\sT$.
Then,  $e=\exp(\Gamma_{\sE}/\Gamma_{\sT})$,
and $\SK(\sE)\cong \mu_s(\sT_0)/\mu_e(\sT_0)$,
where $s = \ind(\sE)$. Moreover,   
there is a short exact sequence 
\begin{equation}\label{maherqde}
 0 \longrightarrow \mathbb Z/{(n,e)}\mathbb Z 
\longrightarrow \hSK(\M_n(\sE)) \xrightarrow{\ \ov \xi\ } 
\SK(\sE) \longrightarrow 0.
 \end{equation} 
\end{enumerate}

\end{theorem}

\begin{proof}
We use the notation in the proof of 
Prop.~\ref{totalfilt}.

Recall from the proof of Prop.~\ref{totalfilt}
that $\sA_h^* = (\sE^*\I_n)\sA_0^*$.  Since $\sA_0^*$ is a 
normal subgroup of $\sA_h^*$, it follows by using 
the 
commutator identity \eqref{commid} that 
$[\sA_h^*,\sA_h^*]\big / [\sA_h^*,\sA_0^*]$ is 
generated by the images of 
$[c\I_n,c'\I_n]= [c,c']\I_n$ for $c,c' \in \sE^*$. 
Note that 
$$
\textstyle
\det_{\sA_0}([c\I_n,c'\I_n])\,=\,
{[c,c']^n}\, [\sE_0^*, \sE_0^*].
$$ 
Furthermore note that the commutators $[c,c']$ generate 
$[\sE^*,\sE^*]$. 
Since the isomorphism $\xi$ maps 
$\ov {[\sA_h^*,\sA_0^*]}$ to 
$\ov {[\sE_h^*,\sE_0^*]}$ by Prop.~\ref{totalfilt}, 
it therefore maps 
maps $ [\sA_h^*,\sA_h^*]\big/[\sA_0^*, \sA_0^*]$
onto $[\sE^*,\sE^*]^n[\sE^*, \sE^*_0]\big/ 
[\sE_0^*,\sE_0^*]$. 
Hence,
$$
\SK^h(\sA) \, = \, \sA^{(1)}/[\sA_h^*, \sA_h^*]
\,\cong \,\sE^{(1)}/[\sE^*,\sE^*]^n[\sE^*, \sE^*_0],
$$
which yields the exact sequence \eqref{mainthwxse}.

For (i)--(iii), recall from 
\cite[Th.~3.4, Lemma~3.5]{HW1} that  there is a
well-defined epimorphism\break 
${\psi\colon\Lambda = \Gamma_{\sE} /\Gamma_{\sT} \wedge 
\Gamma_{\sE} /\Gamma_{\sT} \to 
[\sE^*,\sE^*]/[\sE_0^*,\sE^*]}$,  given as follows:
For $\gamma,\delta\in \Gamma_{\sE}$, take any nonzero $x_\gamma\in 
\sE_\gamma$ and $x_\delta\in \sE_\delta$. Then, 
$$
\psi\big((\gamma+\Gamma_{\sT})\wedge(\delta+\Gamma_{\sT})
\big) \,=\, 
[x_\gamma,x_\delta] \ \text{mod}\ [\sE_0^*, \sE^*].
$$ 
This $\psi$ induces an epimorphism
$$\Lambda/n\Lambda \,\to\, \big([\sE^*,\sE^*]/
[\sE_0^*,\sE^*]\big)
\big/\big([\sE^*,\sE^*]/[\sE_0^*,\sE^*]\big)^n
\,\cong \,
[\sE^*,\sE^*]\big 
/[\sE^*,\sE^*]^n[\sE^*,\sE_0^*],
$$
which yields (i).  Assertion
(ii)  follows  immediately from (i) 
%and \eqref{mainthwxse} 
since the epimorphism 
$\psi$ shows that the exponent of 
$[\sE^*,\sE^*]/[\sE_0^*,\sE^*]$ divides that of $\Lambda$.
Also, (iii) is immediate from (i) and the exact sequence
\eqref{mainthwxse}, since $\Lambda/n\Lambda$ is trivial 
when $\gcd(n,e) = 1$.

For (iv), let $\sE$  be an unramified graded division 
algebra with center $\sT$, i.e., suppose 
$\Gamma_{\sE}=\Gamma_{\sT}$. Then we have 
$\sE^*=\sE_0^*\sT^*$, so $[\sE^*,\sE^*]=[\sE^*,\sE_0^*]$ 
and it follows immediately  from  (\ref{mainthwxse}) 
that   $\hSK(\M_n(\sE)) \cong \SK(\sE)$, for any 
$n \in \mathbb N$. (Compare this with 
Corollary~\ref{homdeun}).  The isomorphism 
$\SK(\sE)\cong \SK(\sE_0)$ for $\sE$ unramified 
is given in \cite[Cor.~3.6(i)]{HW1}.

For (v), let $\sE$  be a totally ramified graded 
division algebra with center $\sT$, i.e.,  $\sE_0=\sT_0$.
 Then ${[\sE^*, \sE_0^*] = [\sE^*,\sT_0^*] = 1}$.   Also, 
by \cite[Prop.~2.1]{hwcor}, 
$[\sE^*,\sE^*]\cong\mu_{e'}(\sT_0)\cong \mathbb Z/e' 
\mathbb Z$, where $e'$ is the exponent of 
the torsion abelian group $\Gamma_{\sE}/ \Gamma_{\sT}$. 
But since $\sE$ is totally ramified, there is a 
nondegenerate symplectic pairing on $\Gamma_\sE/\Gamma_\sT$
induced by commutators in $\sE$ (see \cite[Prop.~2.1,
Remark~2.2(ii)]{hwcor}). 
Hence, $\Gamma_\sE/\Gamma_\sT \cong H\times H$ for some finite abelian
group $H$, which implies that the exponent $e'$ of 
$\Gamma_\sE/\Gamma_\sT$ coincides with the exponent~
$e$ of $\Lambda$.  With this information, 
exact sequence \eqref{maherqde} follows 
from \eqref{mainthwxse}.  The formula for 
$\SK(\sE)$ was given in \cite[Cor.~3.6(ii)]{HW1}
 \end{proof}

\begin{example}\label{trunshift}
For any positive integers $e>1$ and $s$ with $e\mid s$ and 
$s$ having the same prime factors as $e$, it is easy to 
construct examples of graded division algebras $\sE$ with 
center $\sT$ such that $\sE$ is totally ramified over~$\sT$
with $\exp(\Ga_\sE/\Ga_\sT) = e$ and $\ind(\sE) = s$, and 
$\SK(\sE)  \cong \mu_s/\mu_e$.  For example, $\sT$ could be
an iterated Laurent polynomial ring
over the complex numbers,  $\sT = \mathbb C[X_1,X_1\inv, 
X_2, X_2\inv,
\ldots,  X_k, X_k\inv]$ graded by multidegree in 
$X_1, \ldots, X_k$ (so $\Gamma_\sT = \mathbb Z^k$).
For $k$ sufficiently large, one can take $\sE$ to be a
tensor product of suitable graded symbol algebras over
$\sT$, cf. \cite[Ex.~5.3]{HW2}. By choosing $e$ arbitrarily
and choosing $n$ not relatively prime to $e$, one obtains 
explicit examples where $\SK(\M_n(\sE))\not \cong 
\SK(\sE)$ by Th.~\ref{extseq}(v).
\end{example}

The exact sequence~\eqref{mainthwxse}, along with part (i) 
of Th.~\ref{extseq} shows that  $\hSK(\M_n(\sE))$ is a 
finite abelian group with exponent dividing $n\ind(\sE)$
(since $\SK(\sE)$ is finite abelian with exponent dividing
$\ind(\sE)$ by \cite[\S23, Lemma~2]{draxl}).
 However if we  permit shifting in the grading on 
matrices, we can construct more complicated 
reduced Whitehead groups. In the example below we 
construct a simple graded algebra such that its 
homogenous $\SK$ is not even a torsion group when $\sT_0^*$
is not torsion.

\begin{example}\label{hgdstr}
%We construct examples that in this case $\hSK$ is not Morita invariant.  

Let $\sE$  be a graded division algebra totally ramified
over its center $\sT$, with  grade group 
$\Gamma_{\sE}\subseteq \Gamma$. 
Consider $\sA=\M_n(\sE)(\overline \delta)$, where $n>1$
and 
$\overline \delta=(\delta_1,\dots,\delta_n)=
\big(0,\delta,\dots,(n-1)\delta\big)$, with 
$\delta \in \Gamma$ chosen so that the order~$m$ of 
$\delta+\Gamma_{\sE}$  in $\Gamma/\Gamma_{\sE}$ 
exceeds $3n$. Let $s = \ind(\sE)$. We will show that 
\begin{equation}\label{trformula}
\textstyle
\hSK(\M_n(\sE)(\ov \de))\,\cong\,
\big ((\,\prod\limits_{i=1}^{n-1}\sT_0^*\,)\times \mu_s(\sT_0)\big)\big/ H
\quad \text{ where } H = \{(\omega, \ldots, \omega, 
\omega^{2-n}) \mid \omega \in \mu_e(\sT_0)\} \cong \mu_e.
\end{equation}

Note that since the $\delta_i$ are distinct modulo
$\Gamma_\sE$,  the grading on matrices (\ref{mmkkhh}) 
shows that  
$\sA_0$ consists of all diagonal matrices with entries 
from $\sE_0$. We show further that $\Gamma_\sA^*= 
\Gamma_\sE$.  For, recall that $\Gamma_\sA^*$ is 
a subgroup of $\Gamma$ with $\Gamma_\sE \subseteq
\Gamma_\sA^* \subseteq \Gamma_\sA$. From 
\eqref{GammaS}, we have 
$$
\textstyle
\Gamma_\sA \,=\, \bigcup\limits_{i=1}^n\bigcup\limits
_{j=1}^n (\delta_i -\delta_j) +\Ga_\sE \,
= \, \bigcup\limits_{k= -(n-1)}^{n-1} k\delta + \Ga_\sE.  
$$
If $\Gamma_\sA^*\supsetneqq \Gamma_\sE$, then $\ell \delta \in \Gamma_\sA^*$
for some integer $\ell$ with $1 \le |\ell| \le n-1$.
Take the integer $q$ with $n\le q\ell <n+\ell$.  
For any integer $k$ with $ |k|\le n-1$, we have
$$
1\,\le\, q\ell - k \,<\, 2n+\ell -1 \,<\,3n\,\le\, m.
$$
Hence, $(q\ell - k)\delta \notin \Ga_\sE$; so, 
$(q\ell)\delta + \Ga_\sE \ne k\de +\Ga_\sE$ for 
any $k$ with $|k|\le n-1$.  Hence, $q\ell \delta \notin
\Gamma_\sA$,  But, $q \ell \delta$ lies in the group 
$\Ga_\sA^*$, a contradiction.  Thus, $\Ga_\sA^* = \Ga_\sE$.

The formula for $\Gamma_\sA^*$ implies that 
$\sA^*_h=\sA_0^* (\sE^* \I_n)$. Since 
$\sA_0^* = \prod_{i=1}^n \sE_0^* = 
\prod_{i=1}^n \sT_0^*$,
which is abelian and centralized by~$\sE^*\I_n$,  
it follows that 
$[\sA^*_h,\sA^*_h]=[\sE^*,\sE^*]\I_n$. 
By \cite[Prop.~2.1]{hwcor}, 
$
%{[\sA^*_h,\sA^*_h]\cong
[\sE^*,\sE^*]
= \mu_e(\sT_0) = \mu_e$, where $e$ is the exponent of 
the torsion abelian group $\Gamma_{\sE}/ \Gamma_{\sT}$.
Hence, $[\sA^*_h,\sA^*_h] = \mu_e\I_n$.

By Lemma~\ref{indeg0}, $\sA_h^{(1)} \subseteq \sA^*_0 
\subseteq \M_n(\sT_0)$. Now, 
for any matrix $U = \text{diag}(u_1, \ldots, u_n)$ 
$\sA_0^*$, we have 
$$
\Nrd_{\sA_0}(U)\, = \, u_1 \ldots u_n,
$$
so by Lemma~\ref{nomng},  
$$
\Nrd_{\sA} (U)
\, = \, (u_1\dots u_n)^s,
$$ 
where $s=\ind(\sE)$. It follows that 
\begin{align*}
\sA^{(1)}_h\,&=\,
\big\{\text{diag}
(u_1, \ldots, u_n)\mid 
\text{each
$u_i \in \sT_0^*$ and 
$u_1\ldots u_n\in \mu_s(\sT_0)$}\big\} \\
&\cong\, 
\big\{(u_1, \ldots, u_{n-1}, \omega)
\mid \text{each $u_i \in \sT_0^*$
and $\omega \in \mu_s(\sT_0$}\big\}
\, \cong \,\textstyle (\,\prod\limits_{i=1}^{n-1} 
\sT_0^*\,) 
\times  \mu_s(\sT_0).
\end{align*}

In the isomorphism \[\sA_h^{(1)} \cong  
(\,\prod\limits_{i=1}^{n-1} 
\sT_0^*\,) 
\times  \mu_s(\sT_0),\] for any $\omega\in \mu_s(\sT_0)$,
the matrix $\omega\I_n$ maps to $(\omega, \ldots, \omega, 
\omega^{2-n})$.
This yields formula \eqref{trformula} for 
$\SK^h(\sA) = \sA_h^{(1)}/[\sA_h^*, \sA_h^*]$. 
\end{example}

One natural question still unanswered is whether 
{\it inhomogeneous} $\SK$ is Morita invariant in the 
graded setting, i.e., whether for a graded division 
algebra $\sE$,  we have a natural isomorphism 
$\SK(\M_n(\sE)) \cong \SK(\sE)$, for   
$n\in \mathbb N$. This seems to be a difficult 
question, in particular as there does not seem to be a 
notion of (inhomogeneous) Dieudonn\'{e} determinant, which 
is what furnishes the 
Morita isomorphism for  division algebras. 
A key fact which one uses frequently for invertible 
matrices over fields and division rings is that they are 
diagonizable modulo their elementary subgroups. However, 
the work of Bass, Heller and Swan 
(\cite[Lemma 3.2.21]{Rose}) shows that the 
decomposition of an invertible matrix over the graded 
field $F[X,X^{-1}]$ modulo its elementary subgroup  is 
not necessarily diagonal.

\end{document}